\documentclass[global]{svjour}
\usepackage{amssymb}
\usepackage{amsfonts}
\usepackage{amsmath}

\input{tcilatex}

\begin{document}

\title{Applications of the KKM\ property to coincidence theorems,
equilibrium problems, minimax inequalities and variational relation problems}
\author{Monica Patriche}
\institute{University of Bucharest 
\email{monica.patriche@yahoo.com}%
}
\mail{\\
University of Bucharest, Faculty of Mathematics and Computer Science, 14
Academiei Street, 010014 Bucharest, Romania}
\maketitle

\bigskip 

\noindent \textbf{Abstract.}{\small \ In this paper, we establish
coincidence-like results in the case when the values of the correspondences
are not convex. In order to do this, we define a new type of
correspondences, namely} {\small properly quasi-convex-like. Further, we
apply the obtained theorems to solve equilibrium problems and to establish a
minimax inequality. In the last part of the paper, we study the existence of
solutions for generalized vector variational relation problems. Our analysis
is based on the applications of the KKM principle. We establish existence
theorems involving new hypothesis and we improve the results of some recent
papers.\medskip }

\noindent \textbf{Key Words. }KKM property, coincidence theorem; equilibrium
problem; minimax inequality; variational relation problem.{\small \medskip }

\bigskip

\section{\textbf{Introduction}}

The aim of this paper is twofold: firstly, to establish a Fan type geometric
result and to apply it in order to obtain \ some coincidence-like theorems
for the case when the images of the correspondences are not convex. Further,
new theorems concerning the existence of solutions for equilibrium problems
are provided. This study also aims to investigate whether the class of
minimax inequalities can be extended. In fact, we obtain a new general
minimax inequality of the following type: inf$_{x\in X}$sup$_{y\in
Y}t(x,y)\leq \frac{\sup_{y\in Y}\inf_{z\in Z}q(y,z)}{\inf_{z\in Z}\sup_{x\in
X}p(x,z)}.$ Another recent result, due to the author, regarding minimax
inequalities for dicontinuous correspondences, is \cite{pat2}.

In this first part of the article, the originality consists of introducing a
new type of properly quasi-convex-like correspondences, which proved to play
an important role in our results. The method of proof is based on the well
known KKM property. There exists a large literature containing applications
of the KKM property to coincidence theorems, equilibrium theorems, maximal
element theorems and minimax inequalities. We refer the reader, for
instance, to M. Balaj \cite{balaj2}, Lin, Ansari and Wu \cite{lin2003}, Lin
and Wan \cite{lin2004} or Park \cite{park1995}.

Secondly, the paper explores how the KKM principle can promote new more
theorems which show the existence of solutions for some classes of
variational relation problems. We emphasize that, here, the method of
application of the KKM property is new and provides new hypotheses for our
resarch.

The rest of the article is organized as follows. In Section 2, we introduce
notations and preliminary results. In Section 3, a convex-type property for
correspondences is defined and some examples are given, as well. We use this
type of correspondences to obtain coincidence-like theorems, to solve vector
equilibrium problems and to establish a minimax inequality. In Section 4, we
apply the KKM principle to general vector variational relation problems
involving correspondences. Concluding remarks are presented in Section
5.\medskip\ 

\section{\textbf{Preliminaries and notation}}

Throughout this paper, we shall use the following notations and definitions:

Let $A$ be a subset of a vector space $X$, $2^{A}$ denotes the family of all
subsets of $A$ and co$A$ denotes the convex hull of $A$. If $A$ is a subset
of a topological space $X,$ $\overline{A}$ denotes the closure of $A$ in $X$.

If $T$, $G:X\rightarrow 2^{Y}$ are correspondences, then co$G$ and $G\cap
T:X\rightarrow 2^{Y}$ are correspondences defined by $($co$G)(x):=$co$G(x)$
and $(G\cap T)(x):=G(x)\cap T(x),$ for each $x\in X$, respectively.$\medskip 
$

Given a correspondence $T:X\rightarrow 2^{Y},$ for each $x\in X,$ the set $%
T(x)$ is called the \textit{upper section} of $T$ at $x.$ For each $y\in Y,$
the set $T^{-1}(y):=\{x\in X:y\in T(x)\}$ is called the \textit{lower section%
} of $T$ at $y$. The correspondence $T^{-1}:Y\rightarrow 2^{X},$ defined by $%
T^{-1}(y)=\{x\in X:y\in T(x)\}$ for $y\in Y$, is called the \textit{(lower)
inverse} of $T.$

For $A\subset X,$ let $T(A)=\bigcup\nolimits_{x\in A}T(x).$

If $X$ is a nonempty set and $Y$ is a topological space, the correspondence $%
T:X\rightarrow 2^{Y}$ is said to be \textit{transfer open-valued} \cite{tian}
if for any $(x,y)\in X\times Y$ with $y\in T(x),$ there exists an $x^{\prime
}\in X$ such that $y\in $int$T(x^{\prime }).$ $T$ is said to be \textit{%
transfer closed-valued} \cite{tian} if for any $(x,y)\in X\times Y$ with $%
y\notin T(x),$ there exists an $x^{\prime }\in X$ such that $y\notin 
\overline{T(x^{\prime })}.$ The correspondence $T$ is transfer closed-valued
on $X$ if and only if (\cite{tian2}) $\cap _{x\in X}\overline{T}(x)=\cap
_{x\in X}T(x).\medskip $

Further, we present the following lemma (Propostion 1 in \cite{lin2001}).

\begin{lemma}
Let $Y$ be a nonempty set, $X$ be a topological space and $T:X\rightarrow
2^{Y}$ be a correspondence. The following assertions are equivalent:
\end{lemma}

\textit{a) }$T^{-1}:Y\rightarrow 2X$\textit{\ is transfer open-valued and }$%
T $\textit{\ has nonempty values;}

\textit{b) }$X=\bigcup\nolimits_{y\in Y}$int$T^{-1}(y).\medskip $

\textit{Notation }We will denote by $\Delta _{n-1}$ the standard $(n-1)$%
-dimensional simplex in $\mathbb{R}^{n},$ that is, $\Delta _{n-1}:=\left\{
(\lambda _{1},\lambda _{2},...,\lambda _{n})\in \mathbb{R}^{n}:\overset{n}{%
\underset{i=1}{\tsum }}\lambda _{i}=1\text{ and }\lambda _{i}\geqslant
0,i=1,2,...,n\right\} .\medskip $

Let $C^{\ast }(\Delta _{n-1}):=\{g=(g_{1},g_{2},...,g_{n}):\Delta
_{n-1}\rightarrow \Delta _{n-1}$ where $g_{i}$ is continuous, $g_{i}(1)=1$
and $g_{i}(0)=0$ for each\textit{\ }$i\in \{1,2,...,n\}\}.\medskip $

In \cite{pat}, we introduced\ the concept of a weakly naturally
quasi-concave correspondence.

Let $X$ be a nonempty convex subset of a topological vector space $E$\ and $%
Y $ a nonempty subset of a topological vector space $Z$.\textit{\ }The
correspondence\textit{\ }$T:X\rightarrow 2^{Y}$ is said to be weakly
naturally quasi-concave (WNQ) (\cite{pat})\emph{\ }if, for each $n\in 
\mathbb{N}^{\ast }$ and for each finite set $\{x_{1},x_{2},...,x_{n}\}%
\subset X$, there exists $y_{i}\in T(x_{i})$, $(i=1,2,...,n)$ and $g\in
C^{\ast }(\Delta _{n-1}),$ such that $\overset{n}{\underset{i=1}{\tsum }}%
g_{i}(\lambda _{i})y_{i}\in T(\overset{n}{\underset{i=1}{\tsum }}\lambda
_{i}x_{i})$ for every $(\lambda _{1},\lambda _{2},...,\lambda _{n})\in
\Delta _{n-1}.\medskip $

We proved in \cite{pat} the following fixed point result.

\begin{lemma}
( \cite{pat}) \textit{Let }$Y$\textit{\ be a nonempty subset of a
topological vector space }$E,$\textit{\ and }$K$\textit{\ be a }$(n-1)$%
\textit{- dimensional simplex in }$E$\textit{. Let }$T:K\rightarrow 2^{Y}$%
\textit{\ be an weakly naturally quasi-concave correspondence, and }$%
f:Y\rightarrow K$\textit{\ be a continuous function. Then, there exists }$%
x^{\ast }\in K$\textit{\ such that }$x^{\ast }\in f\circ T(x^{\ast })$%
\textit{.\medskip }
\end{lemma}

Now, we recall the generalized KKM mappings, firstly introduced by Park \cite%
{park1989}.

Let $X$ be a convex subset of a linear space, let $Y$ be a topological space
and $T,G:X\rightarrow 2^{Y}$ be two correspondences. We call $G$ a
generalized KKM mapping w.r.t. $T$ if $T($co$A)\subset G(A)$ for each finite
subset $A$ of $X.$ We say that $T$ has the KKM property if $G$ is a
generalized KKM mapping w.r.t. $T$ and the family $\{\overline{G}(x):x\in
X\} $ has the finite intersection property. We denote $KKM(X,Y)=\{T:X%
\rightarrow 2^{Y}:T$ has the KKM\ property$\}.\medskip $

The following lemma is a particular case of Lemma 2.2 in \cite{lin2003}.

\begin{lemma}
Let $X$ be a topological space, $Y$ be a convex set in a topological vector
space. Let $T\in KKM(Y,X)$ be compact and $G:Y\rightarrow 2^{X}$ be a
generalized KKM map w.r.t $T.$ Then, $\overline{T(Y)}\cap
\bigcap\nolimits_{y\in Y}\overline{G(y)}\neq \emptyset .$
\end{lemma}

Let $X$ be a nonempty convex subset of a topological vector space\textit{\ }$%
E,$ $Z$ be a real topological vector space, $Y$ be a subset of $Z$ and $C$
be a pointed closed convex cone in $Z$ with its interior int$C\neq \emptyset
.$ Let $T:X\rightarrow 2^{Z}$ be a correspondence with nonempty values. $T$
is said to be (in the sense of [\cite{ku}, Definition 3.6]) type-(v) \textit{%
properly }$C-$\textit{quasi-convex on} $X~$\cite{ge}, if for any $%
x_{1},x_{2}\in X$ and $\lambda \in \lbrack 0,1],$ either $T(\lambda
x_{1}+(1-\lambda )x_{2})\subset T(x_{1})-C$ or $T(\lambda x_{1}+(1-\lambda
)x_{2})\subset T(x_{2})-C.$

\section{\textbf{Coincidence theorems and applications}}

There are many existed theorems that provide conditions on how to obtain
coincidence points for "adequate" correspondences, that is, correspondences
which satisfy reasonable assumptions concerning the convexity of the images
and continuity. However, there is much less guidance available on how to
obtain similar results under constraints regarding these assumptions. This
section addresses such a challenge and the possible applications of a new
point of view to some classes of generalized vector equilibrium problems and
minimax inequalities.

\subsection{Coincidence theorems and generalized vector equilibrium problems}

In this subsection, we prove some generalized coincidence theorems for the
case when the images of the correspondences are not convex. We work with new
types of properly quasi-convex-like correspondences. By applying our
results, we obtain new theorems concerning the existence of solutions for
generalized vector equilibrium problems.

Now, we present the first result of this subsection. By using Lemma 2, we
establish Theorem 1.

\begin{theorem}
Let $X$ be a simplex in a topological vector space $E,$ let $Y$ be a
Hausdorff space and $T,G:X\rightarrow 2^{Y}$ be correspondences satisfying:
\end{theorem}

\textit{a) }$T$\textit{\ is weakly naturally quasi-concave and compact;}

\textit{b) for each }$y\in T(X),$\textit{\ }$G^{-1}(y)$\textit{\ is convex;}

\textit{c) }$\overline{T(X)}=\cup _{x\in X}$\textit{Int}$G(x).$

\textit{Then, there exists }$x^{\ast }\in X$\textit{\ such that }$T(x^{\ast
})\cap G(x^{\ast })\neq \emptyset .$

\begin{proof}
Since $\overline{T(X)}$ is compact, assumption c) implies that there exists $%
x_{1},x_{2},...,x_{n}\in X$ such that $\overline{T(X)}\subset \cup
_{i=1}^{n} $Int$G(x_{i}).$ We consider $\{\lambda _{1},\lambda
_{2},...\lambda _{n}\}$ the partition of unity corresponding to $\{$int$%
G(x_{i})\}_{i=1,...,n.}.$ We denote $K=$co$\{x_{1},x_{2},...x_{n}\}\subset X$
and we define $f:\overline{T(X)}\rightarrow K$ by

$f(y)=\sum\nolimits_{i=1}^{n}\lambda _{i}(y)x_{i}$ for each $y\in \overline{%
T(X)}.$

We note that $\lambda _{i}(y)\neq \emptyset $ if only if $y\in $int$G(x_{i})$
or, $x_{i}\in G^{-1}(y).$

The continuity of $f$ is obvious and assumption 2) implies $f(y)\in $co$%
\{x_{i}:x_{i}\in G^{-1}(y)\}\subset G^{-1}(y)$ for each $y\in T(X).$

According to Lemma 2, there exists $x^{\ast }\in X$ such that $x^{\ast }\in $
$fT(x^{\ast }).$ In addition, $f^{-1}(x^{\ast })\subset G(x^{\ast })$ and we
obtain that $T(x^{\ast })\cap G(x^{\ast })\neq \emptyset .$
\end{proof}

We define the following type of correspondences.

\begin{definition}
Let $X$ be a nonempty convex subset of a topological vector space\textit{\ }$%
E$ and $Y$ be a real topological vector space$.$ Let $T,G:X\rightarrow 2^{Y}$
be correspondences with nonempty values. $T$ is said to be \textit{properly
quasi-convex }w.r.t. $G$ \textit{on} $X$, if for each $n\in \mathbb{N},$ $%
n\geq 2,$ $x_{1},x_{2},...,x_{n}\subset Y,$ $x\in $co$%
\{x_{1},x_{2},...,x_{n}\},$ there exists $i_{0}\in \{1,2,...,n\}$ such that $%
T(x)\subset G(x_{i_{0}}).$
\end{definition}

\begin{example}
Let $S_{+}^{\prime }((0,0),x):=\{(u,v)\in \lbrack -1,1]\times \lbrack
0,1]:u^{2}+v^{2}\leq x^{2}\}$ and
\end{example}

$S_{-}^{\prime }((0,0),x):=\{(u,v)\in \lbrack -1,1]\times \lbrack
-1,0]:u^{2}+v^{2}\leq x^{2}\}.$

Let us define $T,G:[0,1]\rightarrow 2^{[-1,1]\times \lbrack -1,1]}$ by

$T(x):=\left\{ 
\begin{array}{c}
S_{+}^{\prime }((0,0),x),\text{\ \ \ \ if \ \ \ \ \ }x\in \lbrack 0,1],\text{
}x\neq 1/4; \\ 
S_{-}^{\prime }((0,0),x),\text{ \ \ \ \ \ \ \ \ \ \ \ \ if \ \ \ \ \ \ \ \ \
\ \ \ }x=1/4;%
\end{array}%
\right. $ and

$G(x):=\left\{ 
\begin{array}{c}
S_{+}^{\prime }((0,0),x)\cup \{(x,x):x\in \lbrack 0,1]\},\text{ if \ }x\in
\lbrack 0,1],\text{ }x\neq 1/4; \\ 
S_{-}^{\prime }((0,0),x)\cup \{(x,x):x\in \lbrack -1,0]\},\text{ \ \ \ \ \
if \ \ \ \ \ \ \ \ \ \ \ \ \ }x=1/4.%
\end{array}%
\right. $

Then, $T$ is properly quasi-convex w.r.t. $G.$

\begin{remark}
If $T$ is properly quasi-convex and $T(x)\subseteq G(x)$ for each $x\in X,$
then $T$ is properly quasi-convex\textit{\ }w.r.t. $G.$
\end{remark}

\begin{example}
Let us define $T,G:\mathbb{R}\rightarrow 2^{\mathbb{R}}$ by
\end{example}

$T(x):=\left\{ 
\begin{array}{c}
(-\infty ,4),\text{\ if \ }x\in (-\infty ,2]; \\ 
\lbrack x,3)\text{\ \ if \ \ \ \ \ \ }x\in (2,3); \\ 
(2,\infty ),\text{ if \ }x\in \lbrack 3,\infty );%
\end{array}%
\right. $ and

$G(x):=\left\{ 
\begin{array}{c}
(-\infty ,5),\text{ if \ }x\in (-\infty ,2]; \\ 
\lbrack x,5]\text{ \ \ \ \ if \ \ }x\in (2,3); \\ 
(4,\infty ),\text{ \ if \ }x\in \lbrack 3,\infty ).%
\end{array}%
\right. $

$T$ is properly quasi-convex and $T(x)\subseteq G(x)$ for each $x\in \mathbb{%
R}.$ Then, $T$ is properly quasi-convex w.r.t. $G.$

Example 3 shows that it is not necessary that all images of the
correspondences $T$ to be included in the images of the correspondence $G.$

\begin{example}
Let us define $T,G:[2,3]\rightarrow 2^{\mathbb{R}}$ by
\end{example}

$T(x):=\left\{ 
\begin{array}{c}
(3,\infty ),\text{\ if \ }x=2; \\ 
\lbrack x,3)\text{\ \ if \ \ \ \ \ \ }x\in (2,3); \\ 
(-\infty ,5),\text{ if \ }x=3;%
\end{array}%
\right. $ and

$G(x):=\left\{ 
\begin{array}{c}
(-\infty ,4),\text{ if \ }x=2; \\ 
\lbrack x,3]\text{ \ \ \ \ if \ \ }x\in (2,3); \\ 
(2,\infty ),\text{ \ if \ }x=3.%
\end{array}%
\right. $

Then, $T$ is properly quasi-convex w.r.t. $G.$

\begin{remark}
If $T:X\rightarrow 2^{Y}$ is properly quasi-convex\textit{\ }w.r.t. $%
G:X\rightarrow 2^{Y},$ then, $G$ is a generalized KKM map w.r.t $T.$
\end{remark}

In order to show this assertion, we consider $\{x_{1},x_{2},...,x_{n}\}%
\subset X.$ We want to prove that $T(x)\subset \cup _{i=1}^{n}G(x_{i})$ for
each $x\in $co$\{x_{1},x_{2},...,x_{n}\}.$ Suppose, to the contrary, that
there exist $y^{\ast }\in $co$\{x_{1},x_{2},...,x_{n}\}$ and $y^{\ast }\in
T(x^{\ast })$ such that $y^{\ast }\notin \cup _{i=1}^{n}G(x_{i}).$ Then, $%
y^{\ast }\notin G(x_{i})$ for each $i\in \{1,2,...,n\}$ and $T$ is not
properly quasi-convex\textit{\ }w.r.t. $G.\medskip $

The $T$-properly quasi-convex sets are introduced below.

\begin{definition}
Let $X$ be a topological space, $Y$ be a convex set in a topological vector
space. Let $T:Y\rightarrow 2^{X}$ and $A\subseteq X\times Y.$ The set $A$ is
said to be $T$\textit{-properly quasi-convex} if for each $n\in \mathbb{N},$ 
$n\geq 2,$ $y_{1},y_{2},...,y_{n}\in Y,$ $y\in $co$\{y_{1},y_{2},...,y_{n}\}$
and $x\in T(y),$ there exists $i_{0}\in \{1,2,...,n\}$ such that $%
(x,y_{i_{0}})\in A.$
\end{definition}

\begin{example}
Let us define $T:[0,2]\rightarrow 2^{[-1,1]}$ by
\end{example}

$T(y):=\left\{ 
\begin{array}{c}
1,\text{\ \ \ \ if \ \ \ \ \ }x\in \lbrack 0,2],\text{ }x\neq 1; \\ 
-1,\text{ \ \ \ \ \ \ \ \ \ \ \ \ if \ \ \ \ \ \ \ \ \ \ \ \ }x=1%
\end{array}%
\right. $ and

$A:=\{(0,v):y\in \lbrack 0,1)\cup (1,2)\}\cup \{(x,x):x\in \lbrack
0,1)\}\cup \{(x,x-2):x\in \lbrack 0,1)\}\cup \{(-1,1)\}.$

$A$ is $T$\textit{-}properly quasi-convex.

We emphasize the relation between $T$-properly quasi-convex sets and $T$%
-properly quasi-convex correspondences.

\begin{remark}
We note that if $A$ is $T$-properly quasi-convex and if we define $%
G:X\rightarrow 2^{Y}$ by $G(x)=\{y\in Y:(x,y)\in A\}$ for each $x\in X,$
then, $T$ is properly quasi-convex\textit{\ }w.r.t. $G^{-1}.$
\end{remark}

\begin{example}
In the above example, $G:X\rightarrow 2^{Y}$ is defined by
\end{example}

$G(x)=\{y\in Y:(x,y)\in A\}=\left\{ 
\begin{array}{c}
\lbrack 0,1)\cup (1,2)\text{ \ \ \ \ \ \ if \ \ \ \ \ \ \ \ \ \ }x=1; \\ 
\{(x,x)\}\cup \{(x,x-2)\}\text{ if }x\in \lbrack 0,1); \\ 
\{1\}\text{ \ \ \ \ \ \ \ \ \ \ \ \ if \ \ \ \ \ \ \ \ \ \ \ \ \ \ \ }x=-1%
\end{array}%
\right. $ and

$G^{-1}:Y\rightarrow 2^{X}$ is defined by

$G^{-1}(y)=\{x\in Y:(x,y)\in A\}=\left\{ 
\begin{array}{c}
\{1\}\cup \{y\}\text{ \ \ if \ \ \ }y\in \lbrack 0,1); \\ 
\{1\}\cup \{2-y\}\text{ if }y\in (1,2]; \\ 
\{-1\}\text{ \ \ \ \ \ \ \ \ if \ \ \ \ \ \ \ \ }y=1.%
\end{array}%
\right. $

Then, $T$ is properly quasi-convex\textit{\ }w.r.t. $G^{-1}.\medskip $

A generalization of $T$-properly quasi-convex\textit{\ }sets is introduced
now.

\begin{definition}
Let $X$ be a topological space, $Y$ be a convex set in a topological vector
space. Let $T:Y\rightarrow 2^{X},$ $Q:Y\rightarrow 2^{Y}$ and $A\subseteq
X\times Y.$ The set $A$ is said to be $T$\textit{-properly quasi-convex}
with respect to $Q$ if, for each $n\in \mathbb{N},$ $n\geq 2,$ $%
y_{1},y_{2},...,y_{n}\in Y,$ $y\in $co$\{y_{1},y_{2},...,y_{n}\}$ and $x\in
T(y),$ there exists $i_{0}\in \{1,2,...,n\}$ such that $(x,z)\in A$ for each 
$z\in Q(y_{i_{0}}).$
\end{definition}

We note that if $A$ is $T$-properly quasi-convex with respect to $Q$ and if
we define $G:X\rightarrow 2^{Y}$ by $G(x)=\{y\in Y:P(x)\cap Q(y)=\emptyset
\},$ where $P(x)=\{y\in Y:(x,y)\notin A\}$ for each $x\in X,$ then, $T$ is
properly quasi-convex\textit{\ }w.r.t. $G^{-1}.\medskip $

Theorem 2 is a Fan type geometric result involving $T$-properly quasi-convex
sets.

\begin{theorem}
Let $X$ be a topological space, $Y$ be a convex set in a topological vector
space and let $T\in KKM(Y,X)$ be compact. Let $A\subseteq X\times Y$
satisfying the following conditions: \newline
a) $A$ is $T$\textit{-properly quasi-convex}$;$\newline
b) the correspondence $G:X\rightarrow 2^{Y},$ defined by $G(x)=\{y\in
Y:(x,y)\in A\}$ for each $x\in X,$ is such that $G^{-1}$ is transfer
closed-valued.\newline
Then, there exists $x^{\ast }\in \overline{T(Y)}$ such that $(x^{\ast
},y)\in A$ for all $y\in Y.$
\end{theorem}

\begin{proof}
According to Remarks 2 and 3, $G^{-1}$ is a generalized KKM map w.r.t $T.$
Hence, we can apply Lemma 3 and we obtain that $\overline{T(Y)}\cap
\bigcap\nolimits_{y\in Y}\overline{G^{-1}(y)}\neq \emptyset .$ Consequently,
there exists $x^{\ast }\in \overline{T(Y)}$ such that $(x^{\ast },y)\in A$
for all $y\in Y.$
\end{proof}

We give an example to illustrate the usage of Theorem 2.

\begin{example}
Let us define $T:[0,2]\rightarrow 2^{[-1,1]}$ by
\end{example}

$T(y):=\left\{ 
\begin{array}{c}
1,\text{\ \ \ \ if \ \ \ \ \ }x\in \lbrack 0,2],\text{ }x\neq 1; \\ 
-1,\text{ \ \ \ \ \ \ \ \ \ \ \ \ if \ \ \ \ \ \ \ \ \ \ \ \ }x=1%
\end{array}%
\right. $ and

$A:=\{(0,v):y\in \lbrack 0,2]\}\cup \{(x,x):x\in \lbrack 0,1)\}\cup
\{(x,x-2):x\in \lbrack 0,1)\}\cup \{(-1,1)\}.$

$A$ is $T$\textit{-}properly quasi-convex.

$G^{-1}:[0,2]\rightarrow 2^{[-1,1]}$ is defined by

$G^{-1}(y)=\{x\in \lbrack -1,1]:(x,y)\in A\}=\left\{ 
\begin{array}{c}
\{1\}\cup \{y\}\text{ \ \ if \ \ \ }y\in \lbrack 0,1); \\ 
\{1\}\cup \{2-y\}\text{ if }y\in (1,2]; \\ 
\{-1\}\cup \{1\}\text{ \ \ \ \ if \ \ \ }y=1.%
\end{array}%
\right. $

$G^{-1}$ is transfer closed-valued.

Then, there exists $x^{\ast }=1\in \lbrack -1,1]$ such that $(x^{\ast
},y)\in A$ for all $y\in \lbrack 0,2].\medskip $

Theorem 3 is a Fan type geometric result involving $T$-properly quasi-convex
sets with respect to $Q$.

\begin{theorem}
Let $X$ be a topological space, $Y$ be a convex set in a topological vector
space and let $T\in KKM(Y,X)$ be compact. Let $A\subseteq X\times Y$
satisfying the following conditions: \newline
a) $A$ is $T$\textit{-properly quasi-convex with respect to }$Q:Y\rightarrow
2^{Y};$\newline
b) the correspondence $G:X\rightarrow 2^{Y},$ defined by $G(x)=\{y\in
X:P(x)\cap Q(y)=\emptyset \}$ for each $x\in X,$ where $P(x)=\{y\in
Y:(x,y)\notin A\}$ for each $x\in X,$ is such that $G^{-1}$ is transfer
closed-valued.\newline
Then, there exists $x^{\ast }\in \overline{T(Y)}$ such that $(x^{\ast
},z)\in A$ for all $z\in Q(Y).$
\end{theorem}

\begin{proof}
According to Remarks 2 and 3, $G^{-1}$ is a generalized KKM map w.r.t $T.$
Hence, we can apply Lemma 3 and we obtain that $\overline{T(Y)}\cap
\bigcap\nolimits_{y\in Y}\overline{G^{-1}(y)}\neq \emptyset .$ Consequently,
there exists $x^{\ast }\in \overline{T(Y)}$ such that $(x^{\ast },z)\in A$
for all $z\in Q(Y).$
\end{proof}

Now, we obtain a coincidence-like result, for the case when the convexity of
the images of the correspondence $P$ is missing.

\begin{theorem}
Let $X$ be a topological space and $Y$ be a convex set in a topological
vector space. Let $T\in KKM(Y,X)$ be compact. Let $P:X\rightarrow 2^{Y}$ be
a correspondence such that $X=\cup _{y\in Y}$Int$P^{-1}(y).$ Then, there
exist $n\in \mathbb{N},$ $n\geq 2,$ $B=\{y_{1}^{\ast },y_{2}^{\ast
},...,y_{n}^{\ast }\}\subset Y,$ $y^{\ast }\in $co$\{y_{1}^{\ast
},y_{2}^{\ast },...,y_{n}^{\ast }\}$ and $x^{\ast }\in T(y^{\ast }),$ such
that $y_{i}^{\ast }\in P(x^{\ast })$ for each $i\in \{1,2,...,n\}.$
\end{theorem}

\begin{proof}
Suppose, to the contrary, that for each $n\in \mathbb{N},$ $n\geq 2,$ $%
y_{1},y_{2},...,y_{n}\subset Y,$ $y\in $co$\{y_{1},y_{2},...,y_{n}\}$ and $%
x\in T(y),$ there exists $i_{0}\in \{1,2,...,n\}$ such that $y_{i_{0}}\notin
P(x)$. Then, $(x,y_{i_{0}})\in A,$ where $A=\{(x,y)\in X\times Y:(x,y)\notin 
$Gr$P\}\subseteq X\times Y.$ If we define $G:X\rightarrow 2^{Y}$ by $%
G(x)=\{y\in Y:(x,y)\in A\}=\{y\in Y:y\notin P(x)\},$ we can prove that $%
G^{-1}$ is transfer closed-valued. In order to do this, we notice that the
assumption $X=\cup _{y\in Y}$Int$P^{-1}(y)$ and Lemma 1 imply that $P^{-1}$
is transfer open-valued and $P$ is nonempty valued. The relation between the
correspondences $P$ and $G$ leads us to the conclusion that $G^{-1}$ is
transfer closed-valued. By applying Theorem 2, we obtain that there exists $%
x^{\ast }\in \overline{T(Y)}$ such that $(x^{\ast },y)\in A$ for all $y\in
Y. $ Consequently, $P(x^{\ast })=\emptyset ,$ which is a contradiction.
Hence, there exist $n\in \mathbb{N},$ $n\geq 2,$ $\{y_{1}^{\ast
},y_{2}^{\ast },...,y_{n}^{\ast }\}\subset Y,$ $y^{\ast }\in $co$%
\{y_{1}^{\ast },y_{2}^{\ast },...,y_{n}^{\ast }\}$ and $x^{\ast }\in
T(y^{\ast }),$ such that $y_{i}^{\ast }\in P(x^{\ast })$ for each $i\in
\{1,2,...,n\}.$
\end{proof}

\begin{remark}
If $P(x^{\ast })$ is convex, then, $y^{\ast }\in P(x^{\ast })$ and we obtain
a coincidence theorem.\medskip
\end{remark}

We are establishing the following coincidence-like theorem.

\begin{theorem}
Let $X$ be a topological space and $Y$ be a convex set in a topological
vector space. Let $T\in KKM(Y,X)$ be compact. Let $P:X\rightarrow 2^{Y}$ be
a nonempty valued, lower semicontinuous correspondence and let $%
Q:X\rightarrow 2^{Y}$ be open valud, such that, for each $x\in X,$ there
exists $y\in Y$ such that $P(x)\cap Q(y)\neq \emptyset .$ Then, there exist $%
n\in \mathbb{N},$ $n\geq 2,$ $B=\{y_{1}^{\ast },y_{2}^{\ast
},...,y_{n}^{\ast }\}\subset Y,$ $y^{\ast }\in $co$\{y_{1}^{\ast
},y_{2}^{\ast },...,y_{n}^{\ast }\}$ and $x^{\ast }\in T(y^{\ast }),$ such
that $Q(y_{i}^{\ast })\cap P(x^{\ast })\neq \emptyset $ for each $i\in
\{1,2,...,n\}.$
\end{theorem}

\begin{proof}
Suppose, to the contrary, that for each $n\in \mathbb{N},$ $n\geq 2,$ $%
y_{1},y_{2},...,y_{n}\subset Y,$ $y\in $co$\{y_{1},y_{2},...,y_{n}\}$ and $%
x\in T(y),$ there exists $i_{0}\in \{1,2,...,n\}$ such that for each $z\in
Q(y_{i_{0}}),$ $z\notin P(x)$. Then, for each $n\in \mathbb{N},$ $n\geq 2,$ $%
y_{1},y_{2},...,y_{n}\in Y,$ $y\in $co$\{y_{1},y_{2},...,y_{n}\}$ and $x\in
T(y),$ there exists $i_{0}\in \{1,2,...,n\}$ such that $(x,z)\in A$ for each 
$z\in Q(y_{i_{0}}),$ where $A=\{(x,y)\in X\times Y:(x,y)\notin $Gr$%
P\}\subseteq X\times Y.$

If we define $G:X\rightarrow 2^{Y}$ by $G(x)=\{y\in Y:P(x)\cap
Q(y)=\emptyset \}$ for each $x\in X,$ we can prove that $G^{-1}$ is transfer
closed-valued. $Y\backslash G^{-1}(y)=\{x\in X:P(x)\cap Q(y)\neq \emptyset
\} $ is open. Then, $G^{-1}(y)$ is closed. By applying Theorem 3, we obtain
that there exists $x^{\ast }\in \overline{T(Y)}$ such that $(x^{\ast },y)\in
A$ for all $y\in Q(Y).$ Consequently, $P(x^{\ast })\cap Q(Y)=\emptyset ,$
which is a contradiction. Hence, there exist $n\in \mathbb{N},$ $n\geq 2,$ $%
\{y_{1}^{\ast },y_{2}^{\ast },...,y_{n}^{\ast }\}\subset Y,$ $y^{\ast }\in $%
co$\{y_{1}^{\ast },y_{2}^{\ast },...,y_{n}^{\ast }\}$ and $x^{\ast }\in
T(y^{\ast }),$ such that $Q(y_{i}^{\ast })\cap P(x^{\ast })\neq \emptyset $
for each $i\in \{1,2,...,n\}.$
\end{proof}

Theorem 4 can be generalized in the following way.

\begin{theorem}
Let $X$ be a topological space and $Y$ be a convex set in a topological
vector space. Let $T\in KKM(Y,X)$ be compact. Let $P:X\rightarrow 2^{Y}$ and 
$Q:X\rightarrow 2^{Y}$ be correspondences such that $P$ is nonempty valued
and $X=\cup _{y\in Y}$Int$(G^{\prime })^{-1}(y),$ where $G^{\prime
}:X\rightarrow 2^{Y}$ is defined by $G^{\prime }(x)=\{y\in Y:P(x)\cap
Q(y)\neq \emptyset \}$ for each $x\in X.$ Then, there exist $n\in \mathbb{N}%
, $ $n\geq 2,$ $B=\{y_{1}^{\ast },y_{2}^{\ast },...,y_{n}^{\ast }\}\subset
Y, $ $y^{\ast }\in $co$\{y_{1}^{\ast },y_{2}^{\ast },...,y_{n}^{\ast }\}$
and $x^{\ast }\in T(y^{\ast }),$ such that $Q(y_{i}^{\ast })\cap P(x^{\ast
})\neq \emptyset $ for each $i\in \{1,2,...,n\}.$
\end{theorem}

\begin{proof}
Suppose, to the contrary, that for each $n\in \mathbb{N},$ $n\geq 2,$ $%
y_{1},y_{2},...,y_{n}\subset Y,$ $y\in $co$\{y_{1},y_{2},...,y_{n}\}$ and $%
x\in T(y),$ there exists $i_{0}\in \{1,2,...,n\}$ such that for each $z\in
Q(y_{i_{0}}),$ $z\notin P(x)$. Then, for each $n\in \mathbb{N},$ $n\geq 2,$ $%
y_{1},y_{2},...,y_{n}\in Y,$ $y\in $co$\{y_{1},y_{2},...,y_{n}\}$ and $x\in
T(y),$ there exists $i_{0}\in \{1,2,...,n\}$ such that $(x,z)\in A$ for each 
$z\in Q(y_{i_{0}}),$ where $A=\{(x,y)\in X\times Y:(x,y)\notin $Gr$%
P\}\subseteq X\times Y.$

If we define $G:X\rightarrow 2^{Y}$ by $G(x)=\{y\in Y:P(x)\cap
Q(y)=\emptyset \},$ we can prove that $G^{-1}$ is transfer closed-valued. In
order to do this, we notice that the assumption $X=\cup _{y\in Y}$Int$%
G^{\prime -1}(y)$ and Lemma 1 imply that $(G^{\prime })^{-1}$ is transfer
open-valued. The relation between the correspondences $G^{\prime }$ and $G$
leads us to the conclusion that $G^{-1}$ is transfer closed-valued. By
applying Theorem 2, we obtain that there exists $x^{\ast }\in \overline{T(Y)}
$ such that $(x^{\ast },y)\in A$ for all $y\in Y.$ Consequently, $P(x^{\ast
})=\emptyset ,$ which is a contradiction. Hence, there exist $n\in \mathbb{N}%
,$ $n\geq 2,$ $\{y_{1}^{\ast },y_{2}^{\ast },...,y_{n}^{\ast }\}\subset Y,$ $%
y^{\ast }\in $co$\{y_{1}^{\ast },y_{2}^{\ast },...,y_{n}^{\ast }\}$ and $%
x^{\ast }\in T(y^{\ast }),$ such that $Q(y_{i}^{\ast })\cap P(x^{\ast })\neq
\emptyset $ for each $i\in \{1,2,...,n\}.$
\end{proof}

An application of Theorem 4 is provided in order to establish an existence
result for solutions of a generalized vector equilibrium problem.

We consider the following generalized vector equilibrium problem:

Let $X$ be a topological space, $Y$ be a convex set in a topological vector
space and let $Z$ be a topological vector space. We consider a
correspondence $C:X\rightarrow 2^{Z}$ such that, for each $x\in X,$ $C(x)$
is a closed and convex cone with int$C(x)\neq \emptyset .$ Let $F:X\times
Y\rightarrow 2^{Z}\backslash \{\emptyset \}.$

Find $x^{\ast }\in X$ such that $F(x^{\ast },y)\nsubseteq -$int$C(x^{\ast })$
for each $y\in Y.$

The correspondence $P$ will be needed in our proof. Let $P:X\rightarrow 2^{Y%
\text{ }}$ be defined by $P(x)=\{y\in Y:F(x,y)\subseteq -$int$C(x)\}$ for
each $x\in X$.

We are ready to prove Theorem 7.

\begin{theorem}
Let $X$ be a topological space, $Y$ be a convex set in a topological vector
space and let $Z$ be a topological vector space. Let $T\in KKM(Y,X)$\textit{%
\ such that, for each compact subset }$A$ of $Y,$ $\overline{T(A)}$\textit{\
is compact. }Let $F:X\times Y\rightarrow 2^{Z}\backslash \{\emptyset \}$ and 
$C:X\rightarrow 2^{Z}$ such that, for each $x\in X,$ $C(x)$ is a pointed,
closed and convex cone with int$C(x)\neq \emptyset .$ Assume that:
\end{theorem}

\textit{a) }$T$\textit{\ is properly quasi-convex w.r.t. }$G,$\textit{\
where }$G:X\rightarrow 2^{Y\text{ }}$\textit{\ is defined by }$G(x)=\{y\in
Y:F(x,y)\nsubseteq -$\textit{int}$C(x)\}$\textit{\ for each }$x\in X;$

\textit{b) }$X=\cup _{y\in Y}$\textit{Int}$P^{-1}(y).$

Then, there exists $x^{\ast }\in X$ such that $F(x^{\ast },y)\nsubseteq -$int%
$C(x^{\ast })$ for each $y\in Y.$

\begin{proof}
Suppose, to the contrary, that the considered equilibrium problem does not
have any solutions. It follows that the correspondence $P$ has nonempty
values. According to Theorem 4, there exist $n\in \mathbb{N},$ $n\geq 2,$ $%
\{y_{1}^{\ast },y_{2}^{\ast },...,y_{n}^{\ast }\}\subset Y,$ $y^{\ast }\in $%
co$\{y_{1}^{\ast },y_{2}^{\ast },...,y_{n}^{\ast }\}$ and $x^{\ast }\in
T(y^{\ast }),$ such that $y_{i}^{\ast }\in P(x^{\ast })$ for each $i\in
\{1,2,...,n\}.$ Then, $y_{i}^{\ast }\notin G(x^{\ast })$ for each $i\in
\{1,2,...,n\}$ and, consequently, $T$ is not properly quasi-convex w.r.t. $%
G. $ This fact contradicts a).\medskip
\end{proof}

\begin{theorem}
Let $X$ be a topological space, $Y$ be a convex set in a topological vector
space and let $Z$ be a topological vector space. Let $T\in KKM(Y,X)$\textit{%
\ such that, for each compact subset }$A$ of $Y,$ $\overline{T(A)}$\textit{\
is compact. }Let $F:X\times Y\rightarrow 2^{Z}\backslash \{\emptyset \}$ and 
$C:X\rightarrow 2^{Z}$ such that, for each $x\in X,$ $C(x)$ is a pointed,
closed and convex cone with int$C(x)\neq \emptyset .$ Assume that:
\end{theorem}

\textit{a) }$T$\textit{\ is properly quasi-convex w.r.t. }$G,$\textit{\
where }$G:X\rightarrow 2^{Y\text{ }}$\textit{\ is defined by }$G(x)=\{y\in
Y:F(x,y)\nsubseteq -$\textit{int}$C(x)\}$\textit{\ for each }$x\in X;$

\textit{b) \ for each }$x\in X,$ $F(\cdot ,y):x\rightarrow 2^{Z}\backslash
\{\emptyset \}$ \textit{is u.s.c. with nonempty compact values and the map }$%
W:X\rightarrow 2^{Z}$\textit{\ defined by }$W(x)=Z\backslash ($\textit{-int}$%
C(x))$\textit{\ is u.s.c;}

\textit{c) for each }$x\in X,$\textit{\ there exists }$y\in Y$\textit{\ such
that }$F(x,y)\subseteq -$\textit{int}$C(x).$

Then, there exists $x^{\ast }\in X$ such that $F(x^{\ast },y)\nsubseteq -$int%
$C(x^{\ast })$ for each $y\in Y.$

\begin{proof}
We will prove that for each $y\in Y,$ $P^{-1}(y)$ is open. In order to prove
this, we consider $x\in \overline{X\backslash P^{-1}(y)}$ and a net $%
\{x_{\alpha }\}_{\alpha \in \Lambda }$ in $X\backslash P^{-1}(y)$ such that $%
x_{\alpha }\rightarrow x.$ Since $x_{\alpha }\in X\backslash P^{-1}(y)$ for
each $\alpha \in \Lambda ,$ we have that $F(x_{\alpha },y)\nsubseteq -$int$%
C(x_{\alpha }).$ Then, for each for each $\alpha \in \Lambda ,$ there exists 
$z_{\alpha }\in F(x_{\alpha },y)$ such that $z_{\alpha }\in Z\backslash ($%
-int$C(x_{\alpha })).$ Assumption b) implies that $z\in F(x,y)$ and $z\in
Z\backslash (-$int$C(x_{\alpha })),$ that is, $F(x,y)\nsubseteq -$int$C(x).$
Consequently, $x\in X\backslash P^{-1}(y).$ This shows that $X\backslash
P^{-1}(y)$ is closed and $P^{-1}(y)$ is open for each $y\in Y.$

Assumption c) implies that for each $x\in X,$ $P(x)$ is nonempty. According
to Lemma 1, $X=\cup _{y\in Y}$Int$P^{-1}(y)$ and we can use Theorem 7 in
order to obtain the conclusion.
\end{proof}

Now, by applying Lemma 3, we derive the following coincidence-like theorem.

\begin{theorem}
Let $X$ be a topological space, $Y$ be a convex set in a topological vector
space and Z be a nonempty set.
\end{theorem}

\textit{Let }$P:X\rightarrow 2^{Z},$\textit{\ }$Q:Y\rightarrow 2^{Z}$\textit{%
\ and }$T:Y\rightarrow 2^{X\text{ }}$\textit{\ be correspondences satisfying
the following assumptions}$:$

\textit{a) there exists }$y_{0}\in Y$\textit{\ such that }$P(x)\cap
Q(y_{0})\neq \emptyset $\textit{\ for each }$x\in X;$

\textit{b) }$T\in KKM(Y,X)$\textit{\ is compact.}

\textit{Then, there exist }$n\in \mathbb{N},$ $n\geq 2,$\textit{\ }$%
B=\{y_{1}^{\ast },y_{2}^{\ast },...,y_{n}^{\ast }\}\subset Y,$\textit{\ }$%
y^{\ast }\in $co$\{y_{1}^{\ast },y_{2}^{\ast },...,y_{n}^{\ast }\}$\textit{\
and }$x^{\ast }\in T(y^{\ast }),$\textit{\ such that }$P(x^{\ast })\cap
Q(y_{i}^{\ast })\neq \emptyset $\textit{\ for each} $i\in \{1,2,...,n\}.$

\begin{proof}
Suppose, to the contrary, that for each $n\in \mathbb{N},$ $n\geq 2,\
y_{1},y_{2},...,y_{n}\subset Y,$ $y\in $co$\{y_{1},y_{2},...,y_{n}\}$ and $%
x\in T(y),$ there exists $i_{0}\in \{1,2,...,n\}$ such that $%
Q(y_{i_{0}})\cap P(x)=\emptyset .$ Then, for each $y_{1},y_{2},...,y_{n}%
\subset Y,$ $y\in $co$\{y_{1},y_{2},...,y_{n}\},$ there exists $i_{0}\in
\{1,2,...,n\}$ such that $T(y)\subset G(y_{i_{0}}),$ where $G:Y\rightarrow
2^{X\text{ }}$ is defined by $G(y)=\{x\in X:P(x)\cap Q(y)=\emptyset \}$ for
each $y\in Y.$ We note that $T$ is properly quasi-convex w.r.t. $G.$

We claim that $G$ is generalized KKM w.r.t. $T$. In order to prove this, let
us suppose, to the contrary, that there exist $B=\{y_{1}^{\ast },y_{2}^{\ast
},...,y_{n}^{\ast }\}$ and $x^{\ast }\in T($co$B)\backslash
\bigcup\nolimits_{y^{\ast }\in B}G(y^{\ast }).$ Then, there exists $y^{\ast
}\in $co$B$\textit{\ and }$x^{\ast }\in T(y^{\ast }).$ Since $T$ is properly
quasi-convex w.r.t. $G,$ we have that there exists $i_{0}\in \{1,2,...,n\}$
such that $T(y^{\ast })\subset G(y_{i_{0}}^{\ast }),$ that is, $P(x^{\ast
})\cap Q(y_{i_{0}}^{\ast })=\emptyset .$ Further, we have that $x^{\ast
}\notin \bigcup\nolimits_{i=1}^{n}G(y_{i}^{\ast }).$ Therefore, $P(x^{\ast
})\cap Q(y_{i}^{\ast })\neq \emptyset $ for each $i\in \{1,2,...,n\},$ which
contradicts the last assertion. Consequently, $G$ is generalized KKM w.r.t. $%
T$.

According to Lemma 3, there exists $x^{\ast }\in \bigcap\nolimits_{y\in Y}%
\overline{G(y)}.$ Obviously, $x^{\ast }\in \overline{G(y)}$ for each $y\in Y$
and then, for each $y\in Y,$ there exists a neighborhood $V^{y}$ of $x^{\ast
}$ such that $V^{y}\cap G(y)\neq \emptyset .$ It follows that for each $y\in
Y,$ there exists $x^{y}$ such that $x^{y}\in G(y),$ that is, $P(x^{y})\cap
Q(y)=\emptyset ,$ which contradicts assumption a).\medskip

Theorem 10 is a consequence of Theorem 9. It generalizes some results
concerning the existence of the maximal elements.
\end{proof}

\begin{theorem}
Let $X$ be a topological space, $Y$ be a convex set in a topological vector
space and $Z$ be a nonempty set.
\end{theorem}

\textit{Let }$P:X\rightarrow 2^{Z},$\textit{\ }$Q:Y\rightarrow 2^{Z}$\textit{%
\ and }$T:Y\rightarrow 2^{X\text{ }}$\textit{\ be correspondences satisfying
the following assumptions}$:$

\textit{a) }$T$\textit{\ is properly quasi-convex w.r.t. }$G,$\textit{\ where%
} $G:Y\rightarrow 2^{X\text{ }}$ \textit{is defined by} $G(y)=\{x\in
X:P(x)\cap Q(y)=\emptyset \}$ for each $y\in Y.$

\textit{b) }$T\in KKM(Y,X)$\textit{\ is compact.}

\textit{Then, for each }$y\in Y,$\textit{\ there exists }$x_{0}^{y}$\textit{%
\ such that }$P(x_{0}^{y})\cap Q(y)=\emptyset .$

\begin{proof}
Suppose, to the contrary, that the conclusion of the theorem does not hold.
Therefore, there exists $y_{0}\in Y$ such that $P(x)\cap Q(y_{0})\neq
\emptyset $ for each $x\in X$ and then, assumption a) of Theorem 9 is
fulfilled. By applying Theorem 9, we obtain that there exist $n\in \mathbb{N}%
,$ $n\geq 2,$ $B=\{y_{1}^{\ast },y_{2}^{\ast },...,y_{n}^{\ast }\}\subset Y,$
$y\in $co$B$ and $x^{\ast }\in T(y^{\ast }),$ such that $P(x^{\ast })\cap
Q(y_{i}^{\ast })\neq \emptyset $ for each $i\in \{1,2,...,n\},$ which
contradicts a).
\end{proof}

\begin{remark}
If $Q(y)=Z$ for each $y\in Y,$ then Theorem 10 asserts the existence of the
maximal elements of the correspondence $T.$
\end{remark}

Now, we consider the following generalized vector equilibrium problem:

Let $X$ be a topological space, $Y$ be a convex set in a topological vector
space and $Z$ be a topological vector space.

Let $F:X\times Z\rightarrow 2^{Z}\backslash \{\emptyset \},$ $Q:Y\rightarrow
2^{Z}$ and $C:X\rightarrow 2^{Z}$ such that, for each $x\in X,$ $C(x)$ is a
pointed, closed and convex cone with int$C(x)\neq \emptyset .$ Let $%
F:X\times Z\rightarrow 2^{Z}\backslash \{\emptyset \},$ $Q:Y\rightarrow
2^{Z} $ and $C:X\rightarrow 2^{Z}$ such that, for each $x\in X,$ $C(x)$ is a
pointed, closed and convex cone with int$C(x)\neq \emptyset .$

Then, for each $y\in Y,$ find $x^{\ast }\in X$ such that $F(x^{\ast
},z)\nsubseteq -$int$C(x^{\ast })$ for each $z\in Q(y).$

The correspondence $P$ will be needed in our proof. We define $%
P:X\rightarrow 2^{Z\text{ }}$ by $P(x)=\{z\in Z:F(x,z)\subseteq -$int$C(x)\}$
for each $x\in X$.\medskip

Theorem 11 concerns the existence of solutions for the above generalized
vector equilibrium problem.

\begin{theorem}
Let $X$ be a topological space, $Y$ be a convex set in a topological vector
space and $Z$ be a topological vector space. Let $T\in KKM(Y,X)$\textit{\
such that, for each compact subset }$A$ of $Y,$ $\overline{T(A)}$\textit{\
is compact. }Let $F:X\times Z\rightarrow 2^{Z}\backslash \{\emptyset \},$ $%
Q:Y\rightarrow 2^{Z}$ and $C:X\rightarrow 2^{Z}$ such that, for each $x\in
X, $ $C(x)$ is a pointed, closed and convex cone with int$C(x)\neq \emptyset
.$ Assume that $T$ is \textit{properly quasi-convex }w.r.t. $G$, where $%
G:Y\rightarrow 2^{X\text{ }}$ is defined by $G(y)=\{x\in X:P(x)\cap
Q(y)=\emptyset \}$ for each $y\in Y$.
\end{theorem}

\textit{Then, for each }$y\in Y,$\textit{\ there exists }$x^{\ast }\in X$%
\textit{\ such that }$F(x^{\ast },z)\nsubseteq -$\textit{int}$C(x^{\ast })$%
\textit{\ for each} $z\in Q(y).$

\begin{proof}
According to Theorem 10, for each $y\in Y,$\ there exists $x^{\ast }$\ such
that $P(x^{\ast })\cap Q(y)=\emptyset ,$ that is, for each $y\in Y,$\ there
exists $x^{\ast }$\ such that $R(x^{\ast },y)=\{z\in Q(y):F(x^{\ast
},z)\subseteq -$int$C(x^{\ast })\}=\emptyset .$

Consequently, for each $y\in Y,$ there exists $x^{\ast }\in X$ such that $%
F(x^{\ast },z)\nsubseteq -$int$C(x^{\ast })$ for each $z\in Q(y).$
\end{proof}

A new coincidence-like theorem can be proved under new assumptions.

\begin{theorem}
Let $X$ be a topological space, $Y$ be a convex set in a topological vector
space and Z be a nonempty set.
\end{theorem}

\textit{Let }$P:X\rightarrow 2^{Z},$\textit{\ }$Q:Y\rightarrow 2^{Z}$\textit{%
\ and }$T:Y\rightarrow 2^{X\text{ }}$\textit{\ be correspondences satisfying
the following assumptions}$:$

\textit{a) for each }$x\in X,$ \textit{there exists }$y\in Y$\textit{\ such
that }$P(x)\cap Q(y)\neq \emptyset ;$

\textit{b) for each }$x\in X$\textit{\ and }$y\in Y$ such that $P(x)\cap
Q(y)\neq \emptyset ,$\textit{\ there xists }$y^{\prime }\in Y$\textit{\ and
a neighborhood }$U$\textit{\ of }$x$\textit{\ such that }$P(x^{\prime })\cap
Q(y^{\prime })\neq \emptyset $\textit{\ for each }$x^{\prime }\in U;$

\textit{c) }$T\in KKM(Y,X)$\textit{\ is compact.}

\textit{Then, there exist }$n\in \mathbb{N},$ $n\geq 2,$\textit{\ }$%
B=\{y_{1}^{\ast },y_{2}^{\ast },...,y_{n}^{\ast }\}\subset Y,$\textit{\ }$%
y^{\ast }\in $co$B$\textit{\ and }$x^{\ast }\in T(y^{\ast }),$\textit{\ such
that }$P(x^{\ast })\cap Q(y_{i}^{\ast })\neq \emptyset $\textit{\ for each} $%
i\in \{1,2,...,n\}.$

\begin{proof}
Suppose, to the contrary, that for each $n\in \mathbb{N},$ $n\geq 2,\
y_{1},y_{2},...,y_{n}\subset Y,$ $y\in $co$\{y_{1},y_{2},...,y_{n}\}$ and $%
x\in T(y),$ there exists $i_{0}\in \{1,2,...,n\}$ such that $%
Q(y_{i_{0}})\cap P(x)=\emptyset .$ Then, for each $y_{1},y_{2},...,y_{n}%
\subset Y,$ $y\in $co$\{y_{1},y_{2},...,y_{n}\},$ there exists $i_{0}\in
\{1,2,...,n\}$ such that $T(y)\subset G(y_{i_{0}}),$ where $G:Y\rightarrow
2^{X\text{ }}$ is defined by $G(y)=\{x\in X:P(x)\cap Q(y)=\emptyset \}$ for
each $y\in Y.$ We note that $T$ is properly quasi-convex w.r.t. $G.$

We claim that $G$ is generalized KKM w.r.t. $T$. In order to prove this, let
us suppose, to the contrary, that there exist $n\in \mathbb{N},$ $n\geq 2,$ $%
B=\{y_{1}^{\ast },y_{2}^{\ast },...,y_{n}^{\ast }\}$ and $x^{\ast }\in T($co$%
B)\backslash \bigcup\nolimits_{y^{\ast }\in B}G(y^{\ast }).$ Then, there
exists $y^{\ast }\in $co$B$\textit{\ and }$x^{\ast }\in T(y^{\ast }).$ Since 
$T$ is properly quasi-convex w.r.t. $G,$ we have that there exists $i_{0}\in
\{1,2,...,n\}$ such that $T(y^{\ast })\subset G(y_{i_{0}}^{\ast }),$ that
is, $P(x^{\ast })\cap Q(y_{i_{0}}^{\ast })=\emptyset .$ Further, we have
that $x^{\ast }\notin \bigcup\nolimits_{i=1}^{n}G(y_{i}^{\ast }).$
Therefore, $P(x^{\ast })\cap Q(y_{i}^{\ast })\neq \emptyset $ for each $i\in
\{1,2,...,n\},$ which contradicts the last assertion. Consequently, $G$ is
generalized KKM w.r.t. $T$.

According to Lemma 3, there exists $x^{\ast }\in \bigcap\nolimits_{y\in Y}%
\overline{G(y)}.$ Assumption b) implies that $G$ is transfer closed-valued
and then, $x^{\ast }\in \bigcap\nolimits_{y\in Y}G(y).$ It follows that $%
P(x^{\ast })\cap Q(y)=\emptyset $ for each $y\in Y,$ which contradicts a)$.$%
\medskip

Theorem 13 is a consequence of Theorem 12. It generalizes some results
concerning the existence of the maximal elements.
\end{proof}

\begin{theorem}
Let $X$ be a topological space, $Y$ be a convex set in a topological vector
space and $Z$ be a nonempty set.
\end{theorem}

\textit{Let }$P:X\rightarrow 2^{Z},$\textit{\ }$Q:Y\rightarrow 2^{Z}$\textit{%
\ and }$T:Y\rightarrow 2^{X\text{ }}$\textit{\ be correspondences satisfying
the following assumptions}$:$

\textit{a) }$T$\textit{\ is properly quasi-convex w.r.t. }$G,$\textit{\ where%
} $G:Y\rightarrow 2^{X\text{ }}$ \textit{is defined by} $G(y)=\{x\in
X:P(x)\cap Q(y)=\emptyset \}$ for each $y\in Y.$

\textit{b) }$T\in KKM(Y,X)$\textit{\ is compact.}

\textit{Then}$,$\textit{\ there exists }$x_{0}$\textit{\ such that }$%
P(x_{0})\cap Q(y)=\emptyset $ \textit{for each }$y\in Y.$

\begin{proof}
Suppose, to the contrary, that the conclusion of the theorem does not hold.
Therefore, for each $x\in X,$ there exists $y_{0}\in Y$ such that $P(x)\cap
Q(y_{0})\neq \emptyset $ and then, assumption a) of Theorem 12 is fulfilled.
By applying Theorem 12, we obtain that there exist $n\in \mathbb{N},$ $n\geq
2,$ $B=\{y_{1}^{\ast },y_{2}^{\ast },...,y_{n}^{\ast }\}\subset Y,$ $y^{\ast
}\in $co$B$ and $x^{\ast }\in T(y^{\ast }),$ such that $P(x^{\ast })\cap
Q(y_{i}^{\ast })\neq \emptyset $ for each $i\in \{1,2,...,n\},$ which
contradicts a).
\end{proof}

\begin{remark}
If $Q(y)=Z$ for each $y\in Y,$ then Theorem 13 asserts the existence of the
maximal elements of the correspondence $T.$
\end{remark}

\subsection{A \textbf{minimax inequality}}

The classical Ky Fan's minimax inequalities \cite{fan1}-\cite{fan3} have
played an important role in the study of modern nonlinear analysis. This
subsection is devoted to the research of a generalized Ky Fan minimax
inequality for vector-valued functions.

First, we prove the following interesting theorem.

\begin{theorem}
Let $X$ be a topological space, $Y$ and $Z$ be convex sets in topological
spaces, $p:X\times Z\rightarrow \mathbb{R},$ $q:X\times Z\rightarrow \mathbb{%
R},$ $q:X\times Y\rightarrow \mathbb{R}$ functions and $\alpha ,\beta
,\lambda $ real numbers. Suppose that the following conditions are fulfilled:
\end{theorem}

\textit{a) }$\mathit{p}$\textit{\ is transfer upper semicontinuous in }$%
\mathit{x}$\textit{;}

\textit{b) for each }$x\in X,$\textit{\ }$y_{1},y_{2},...,y_{n}\in Y,$%
\textit{\ }$y\in $co$\{y_{1},y_{2},...,y_{n}\}$\textit{\ if }$t(x,y)\geq
\lambda ,$\textit{\ then there exists }$i_{0}\in \{1,2,...,n\}$\textit{\
such that }$t(x,y_{i_{0}})\geq \lambda ;$

\textit{c) for each }$x\in X,$\textit{\ }$y\in Y$\textit{\ and }$z\in Z,$%
\textit{\ if }$t(x,y)\geq \lambda ,$\textit{\ then }$p(x,z)<\alpha $\textit{%
\ and }$q(y,z)>\beta ;$

\textit{d)} \textit{for each }$x\in X$ \textit{and }$z\in Z$\textit{\ such
that }$p(x,z)<\alpha ,$\textit{\ there exists }$y\in Y$ \textit{such that} $%
q(y,z)<\beta ;$

\textit{e) there exists a compact subset }$K$\textit{\ of }$X$\textit{\ such
that for each }$x\in X\backslash K,$\textit{\ }$t(x,y)<\lambda $\textit{\
for all }$y\in Y;$

\textit{f) the map }$T:Y\rightarrow 2^{X},$\textit{\ defined by }$%
T(y)=\{x\in X:t(x,y)\geq \lambda \}$\textit{\ for each }$y\in Y,$ \textit{%
has the KKM property.}

\textit{Then, for each }$y\in Y,$ \textit{there exists }$x_{0}\in X$\textit{%
\ such that for }$z\in Z$\textit{\ with the property that }$q(y,z)<\beta ,$%
\textit{\ it is true that} $p(x_{0},z)\geq \alpha .$

\begin{proof}
We start the proof by defining the correspondences $P:X\rightarrow 2^{Z},$ $%
Q:Y\rightarrow 2^{Z}$ and $T:Y\rightarrow 2^{X}$ by

$P(x)=\{z\in Z:p(x,z)<\alpha \},$ $Q(y)=\{z\in Z:q(y,z)<\beta \}$ and $%
T(y)=\{x\in X:t(x,y)\geq \lambda \}$ for each $x\in X,$ respectively $y\in
Y. $

According to assumption e), $T(Y)\subset K$. It follows that $T$ is a
compact correspondence. In addition, according to assumption d), $T$ has the
KKM property.

We claim that $T$ is properly quasi-convex w.r.t. $G,$ where $G:Y\rightarrow
2^{X\text{ }}$ is defined by $G(y)=\{x\in X:P(x)\cap Q(y)=\emptyset \}$ for
each $y\in Y.$

In order to prove this, let us consider $y_{1},y_{2},...,y_{n}\in Y,$ $y\in $%
co$\{y_{1},y_{2},...,y_{n}\}$ and $x\in T(y).$ Then, $t(x,y)\geq \lambda .$

Assumption b) implies that there exists $i_{0}\in \{1,2,...,n\}$ such that $%
t(x,y_{i_{0}})\geq \lambda .$ Further, assumption c) holds and we conclude
that for each $z\in Z,$ $p(x,z)<\alpha $ and $q(y_{i_{0}},z)>\beta ,$ that
is, $P(x)\cap Q(y_{i_{0}})=\emptyset .$ Therefore, $x\in G(y_{i_{0}}),$
which implies that there exists $i_{0}\in \{1,2,...,n\}$ such that $%
T(x_{\lambda })\subset G(y_{i_{0}}).$ We proved that $T$ is properly
quasi-convex w.r.t. $G.$

Based on assumptions a) and d), we conclude that $G$ is transfer
closed-valued.

All assumptions of Theorem 13 are fulfilled. By applying this result, we
obtain that there exists $x_{0}\in X$ such that $P(x_{0})\cap Q(y)=\emptyset 
$ for each $y\in Y.$ Consequently, there exists $x_{0}\in X$ such that for
each $y\in Y$ and $z^{y}\in Z$ with the property that $q(y,z^{y})<\beta ,$
it is true that $p(x_{0},z^{y})\geq \alpha .$
\end{proof}

Theorem 15 is a consequence of Theorem 14.

\begin{theorem}
Let $X$ be a topological space, $Y$ and $Z$ be convex sets in topological
spaces and let $p:X\times Z\rightarrow \mathbb{R}_{+},$ $q:X\times
Z\rightarrow \mathbb{R}_{+},$ $q:X\times Y\rightarrow \mathbb{R}_{+}$ be
functions. Suppose that the following conditions are fulfilled:
\end{theorem}

\textit{a) }$\mathit{p}$\textit{\ is transfer upper semicontinuous in }$%
\mathit{x}$\textit{;}

\textit{b) for each }$x\in X,$\textit{\ }$y_{1},y_{2},...,y_{n}\in Y,$%
\textit{\ }$y\in $co$\{y_{1},y_{2},...,y_{n}\},$\textit{\ if }$t(x,y)\geq
\lambda ,$\textit{\ then there exists }$i_{0}\in \{1,2,...,n\}$\textit{\
such that }$t(x,y_{i_{0}})\geq \lambda ;$

\textit{c)} \textit{for each }$x\in X$ \textit{and }$z\in Z$\textit{\ such
that }$p(x,z)<\alpha ,$\textit{\ there exists }$y\in Y$ such that $%
q(y,z)<\beta ;$

\textit{d) for each }$x\in X,$\textit{\ }$y\in Y$\textit{\ and }$z\in Z,$%
\textit{\ }$t(x,y)\leq \frac{q(y,z)}{p(x,z)};$

\textit{e) for each }$\lambda <$\textit{inf}$_{x\in X}$\textit{sup}$_{y\in
Y}t(x,y),$\textit{\ the map }$T:Y\rightarrow 2^{X},$\textit{\ defined by }$%
T(y)=\{x\in X:t(x,y)\geq \lambda \}$\textit{\ for each }$y\in Y,$ \textit{%
has the KKM property.}

\textit{Then, inf}$_{x\in X}$\textit{sup}$_{y\in Y}t(x,y)\leq \frac{%
\sup_{y\in Y}\inf_{z\in Z}q(y,z)}{\inf_{z\in Z}\sup_{x\in X}p(x,z)}.$

\begin{proof}
Let us suppose that the conclusion does not hold. Hence, inf$_{x\in X}$sup$%
_{y\in Y}t(x,y)>\frac{\sup_{y\in Y}\inf_{z\in Z}q(y,z)}{\inf_{z\in
Z}\sup_{x\in X}p(x,z)}.$

The constants $\alpha ,\beta ,\lambda \in \mathbb{R}_{+}$ can be chosen such
that inf$_{x\in X}$sup$_{y\in Y}t(x,y)>\lambda ,$ $\sup_{y\in Y}\inf_{z\in
Z}q(y,z)<\beta $, $\inf_{z\in Z}\sup_{x\in X}p(x,z)>\alpha $ and $\lambda >%
\frac{\beta }{\alpha }.$

We claim that condition c) of the above theorem is fulfilled. Indeed, let $%
x\in X,$ $y\in Y$ such that $t(x,y)\geq \lambda >\frac{\beta }{\alpha }.$
According to assumption d) of Theorem 14, $\frac{q(y,z)}{p(x,z)}>\frac{\beta 
}{\alpha }$ for each $z\in Z.$ Then, there exists $\varepsilon >0$ such that 
$p(x,z)<\alpha \varepsilon $ and $q(y,z)>\beta \varepsilon $ for each $z\in
Z.$

All assumptions of Theorem 14 are fulfilled. By applying this result, we
obtain that there exists $x_{0}\in X$ such that for each $y\in Y$ and for $%
z^{y}\in Z$ with the property that $q(y,z^{y})<\beta \varepsilon ,$ it is
true that $p(x_{0},z^{y})\geq \alpha \varepsilon .$

Therefore, there exists $x_{0}\in X$ such that for each $y\in Y$ there
exists $z^{y}$ with the property that $t(x_{0},y)<\frac{q(y,z^{y})}{%
p(x_{0}^{y},z^{y})}<\frac{\beta \varepsilon }{\alpha \varepsilon }=\frac{%
\beta }{\alpha }<\lambda $. It follows that inf$_{x\in X}$sup$_{y\in
Y}t(x,y)\leq \lambda ,$ which contradicts the choice of $\lambda .$
\end{proof}

\section{\textbf{Applications of the KKM principle to general vector
variational inclusion problems involving correspondences }}

The vector variational inclusion problem is considered to be a model which
unifies several other problems, for instance, vector variational
inequalities, vector optimization problems, equilibrium problems or fixed
points theorems. For further information on this topic, the reader is
referred to the following list of selected publications: \cite{ans3}, \cite%
{balaj}-\cite{dinh}, \cite{g}, \cite{kim1987}, \cite{lin}, \cite{peng}, \cite%
{sach}, \cite{zhang}, \cite{zheng}.

We report new results concerning the existence of solutions for the general
vector variational inclusion problems under new assumptions. {\normalsize %
The methodology of the proofs relies on the applications of the KKM
principle. }

In this subsection, we will use a particular form of the KKM principle. We
start by presenting it here. We note that its open version is due to Kim 
\cite{kim1987} and Shih and Tan \cite{shih}.\medskip

Let $X$ be a subset of a topological vector space and $D$ a nonempty subset
of $X$ such that co$D\subset X.$

$T:D\rightarrow 2^{X}$ is called a \textit{KKM correspondence} if co$%
N\subset T(N)$ for each $N\in \langle D\rangle ,$ where $\langle D\rangle $
denotes the class of all nonempty finite subsets of $D$.\medskip

\textbf{KKM principle }Let $D$ be a set of vertices of a simplex $S$ and $%
T:D\rightarrow 2^{S}$ a correspondence with closed (respectively open)
values such that

co$N\subset T(N)$ for each $N\subset D.$

Then, $\tbigcap\nolimits_{z\in D}T(z)\neq \emptyset .\medskip $

As application of the KKM principle, we recall the following lemma.

\begin{lemma}
Let $X$ be a subset of a topological vector space, $D$ a nonempty subset of $%
X$ such that co$D\subset X$ and $T:D\rightarrow 2^{X}$ a \textit{KKM
correspondence with closed (respectively open) values. Then }$\{T(z)\}_{z\in
D}$ has the finite intersection property.
\end{lemma}

\subsection{Main results}

In this subsection, we study the existence of solutions for the following
types of variational relation problems. We emphasize that the method of
application of the KKM property is new.\medskip

Let $X$ be a nonempty subset of a topological vector space $E$ and let $Z$
be a topological vector space. Let $S_{1},S_{2}:X\rightarrow 2^{X},$\textit{%
\ }$T:X\times X\rightarrow 2^{X}$,\textit{\ }$F:T(X\times X)\times X\times
X\rightarrow 2^{Z}$\textit{\ }and\textit{\ }$G:T(X\times X)\times X\times
X\rightarrow 2^{Z}$ be correspondences with nonempty values. We consider the
following generalized vector problems:

(IP 1): Find $x^{\ast }\in S_{1}(x^{\ast })\ $such that\textit{\ }$\forall
y\in S_{2}(x^{\ast }),\forall t^{\ast }\in T(x^{\ast },y),F(t^{\ast
},y,x^{\ast })\subseteq G(t^{\ast },x^{\ast },x^{\ast })$ and

(IP 2): Find $x^{\ast }\in S_{1}(x^{\ast })\ $such that\textit{\ }$\forall
y\in S_{2}(x^{\ast }),\forall t^{\ast }\in T(x^{\ast },y),F(t^{\ast
},y,x^{\ast })\cap G(t^{\ast },x^{\ast },x^{\ast })\neq \emptyset .\mathit{.}%
\medskip $

For a motivation and special cases of these considered problems, the reader
is referred to \cite{hai}. We note that this problem also generalizes the
vector equilibrium problems considered in Subsection 3.1\medskip

Now, we are establishing an existence theorem for the general vector
variational inclusion problem IP 1, by using Lemma 4.

\begin{theorem}
\textit{Let }$Z$\textit{\ be a Hausdorff topological vector space and let }$%
X $\textit{\ be a nonempty compact convex subset of a\ topological vector
space }$E$\textit{. Let }$S_{1},S_{2}:X\rightarrow 2^{X},$\textit{\ }$%
T:X\times X\rightarrow 2^{X}$,\textit{\ }$F:T(X\times X)\times X\times
X\rightarrow 2^{Z}$\textit{\ and }$G:T(X\times X)\times X\times X\rightarrow
2^{Z}$ \textit{be correspondences with nonempty values. Assume that:}
\end{theorem}

\textit{i) }$S_{1}^{-1}$\textit{\ and }$S_{2}^{-1}$\textit{\ are open-valued
and for each }$u\in X,$\textit{\ the set }$\{x\in X:\exists t\in
T(x,u),F(t,u,x)\nsubseteq G(t,x,x)\}$\textit{\ is open; }

\textit{ii)} $S_{1}$\textit{\ and }$S_{2}$\textit{\ are convex valued and
for each }$x\in X,$\textit{\ the set }$\{u\in X:\exists t\in
T(x,u),F(t,u,x)\nsubseteq G(t,x,x)\}$\textit{\ is convex; }

\textit{iii)} \textit{the set }$A=\{x\in X:$\textit{\ there exist }$u\in
S_{2}(x)$ \textit{and }$t\in T(x,u)$\textit{\ such that} $F(t,u,x)\nsubseteq
G(t,x,x)\}$ \textit{is closed;}

\textit{iv) }$\mathit{\{x}\in X:x\in S_{2}(x)\}=\emptyset ;$

\textit{v) there exists }$M\in \langle A\rangle $\textit{\ such that }$%
\bigcup\nolimits_{u\in M}[\{x\in X:\exists t\in T(x,u),F(t,u,x)\nsubseteq
G(t,x,x)\}\cap S_{2}^{-1}(u)=X$ or $\bigcup\nolimits_{u\in M}[(X\backslash
A)\cap \mathit{\ }S_{1}^{-1}(u)]=X.$

\textit{Then, there exists }$x^{\ast }\in S_{1}(x^{\ast })$\ \textit{such
that }$\forall y\in S_{2}(x^{\ast }),$ $\forall t^{\ast }\in T(x^{\ast },y),$
$F(t^{\ast },y,x^{\ast })\subseteq G(t^{\ast },x^{\ast },x^{\ast })$\textit{%
.\medskip }

\textit{Proof.} Let $P,G:X\rightarrow 2^{X}$ be defined by

$P(x)=\{u\in X:\exists t\in T(x,u),F(t,u,x)\nsubseteq G(t,x,x)\},$ for each $%
x\in X$ and

$G(x)=S_{2}(x)\cap P(x),$ for each $x\in X.$

We are going to show that there exists $x^{\ast }\in X$ such that $x^{\ast
}\in S_{1}(x^{\ast })$ and $S_{2}(x^{\ast })\cap P(x^{\ast })=\emptyset .$

We consider two cases.

Case I.

$A=\{x\in X:$\textit{\ }$P(x)\cap S_{2}(x)\neq \emptyset \}=\{x\in X:$%
\textit{\ }$G(x)\neq \emptyset \}$\textit{\ }is nonempty.

The correspondence $G:A\rightarrow 2^{X},$ defined by $G(x)=S_{2}(x)\cap
P(x) $ for each $x\in A,$ is nonempty valued on $A.$ It is obvious that for
each $u\in X,$ $G^{-1}(u)=P^{-1}(u)\cap S_{2}^{-1}(u)$ is a convex set as
intersection of convex sets.

Further, let us define the correspondence $H:X\rightarrow 2^{X}$ by

$H(x)=\left\{ 
\begin{array}{c}
G(x),\text{ if }x\in A; \\ 
S_{1}(x),\text{\ otherwise}%
\end{array}%
\right. =\left\{ 
\begin{array}{c}
S_{2}(x)\cap P(x),\text{ if }x\in A; \\ 
S_{1}(x),\text{ \ \ \ otherwise.}%
\end{array}%
\right. $

According to ii), $H$ is convex valued.

For each $u\in X,$

$H^{-1}(u)=\{x\in X:$\textit{\ }$\mathit{u}\in H(x)\}=$

\ \ \ \ \ \ \ \ \ \ \ $=\{x\in A:$\textit{\ }$\mathit{u}\in S_{2}(x)\cap
P(x)\}\cup \{x\in X\setminus A:$\textit{\ }$\mathit{u}\in S_{1}(x)\}=$

$\ \ \ \ \ \ \ \ =[P^{-1}(u)\cap S_{2}^{-1}(u)]\cup \lbrack (X\backslash
A)\cap $\textit{\ }$S_{1}^{-1}(u)].$

According to i) and iii), $H^{-1}$ is open-valued.

If v) is satisfied, then, there exists $M\in \langle X\rangle $ such that $%
\tbigcup\nolimits_{x\in M}H^{-1}(x)=X.$

Let us define $Q:X\rightarrow 2^{X}$ by $Q(x):=X\backslash H^{-1}(x)$ for
each\textit{\ }$x\in X.$

The correspondence $Q$ is closed-valued and $\tbigcap\nolimits_{x\in
M}Q(x)=X\backslash \tbigcup\nolimits_{x\in M}H^{-1}(x)=\emptyset .$

According to Lemma 4, we conclude that $Q$ is not a KKM correspondence.
Thus, there exists $N\in \langle X\rangle $ such that\textit{\ }co$%
N\varsubsetneq Q(N)=\tbigcup\nolimits_{x\in N}(X\backslash H^{-1}(x)).$

Hence, there exists $x^{\ast }\in $co$N$ with the property that $x^{\ast
}\in H^{-1}(x)$ for each $x\in N,$ which implies $N\subset H(x^{\ast })$. It
is clear that co$N\subset $co$H(x^{\ast })=H(x^{\ast })$. Consequently, $%
x^{\ast }\in $co$N\subset $co$H(x^{\ast })=H(x^{\ast }),$ which means that $%
x^{\ast }\in H(x^{\ast }),$ that is, $x^{\ast }$ is a fixed point for $H.$

We notice that, if $x^{\ast }\in A$, then, $x^{\ast }\in S_{2}(x^{\ast
})\cap P(x^{\ast }),$ which contradicts iv). Therefore, $x^{\ast }\in
X\backslash A$ and $x^{\ast }\in S_{1}(x^{\ast }).$ Since $x^{\ast }\in
X\backslash A,$ we conclude that $S_{2}(x^{\ast })\cap P(x^{\ast
})=\emptyset .$ This shows that $\forall y\in S_{2}(x^{\ast }),$ $\forall
t^{\ast }\in T(x^{\ast },y),$ $F(t^{\ast },y,x^{\ast })\subseteq G(t^{\ast
},x^{\ast },x^{\ast }).$

Case II.

$A=\{x\in X:$\textit{\ }$P(x)\cap S_{2}(x)\neq \emptyset \}=\{x\in X:$%
\textit{\ }$G(x)\neq \emptyset \}\mathit{=}\emptyset $.

In this case, $G(x)=\emptyset $ for each $x\in X.$

Let us define $Q:X\rightarrow 2^{X}$ by $Q(x):=X\backslash S_{1}(x)$ for each%
\textit{\ }$x\in X.$ The proof follows the same line as above and we obtain
that there exists $x^{\ast }\in X$ such that $x^{\ast }\in S_{1}(x^{\ast }).$

Obviously, $G(x^{\ast })=\emptyset .$ Consequently, the conclusion holds in
Case II.

\begin{remark}
Assumption i) can be replaced with
\end{remark}

i') \textit{\ }$S_{1}^{-1}$\ and $S_{2}^{-1}$ are closed-valued and for each 
$u\in X,$\ the set\textit{\ }$\{x\in X:\exists t\in
T(x,u),Q(t,u,x)\nsubseteq G(t,x,x)\}$ is closed.\textit{\ }

In this case, $Q$ is open-valued and the open version of Lemma 4 can be
applied.\smallskip

Using a similar argument as in the proof of the above result, we obtain the
following theorem.

\begin{theorem}
\textit{Let }$Z$\textit{\ be a Hausdorff topological vector space and let }$%
X $\textit{\ be a nonempty compact convex subset of a\ topological vector
space }$E$\textit{. Let }$S_{1},S_{2}:X\rightarrow 2^{X},$\textit{\ }$%
T:X\times X\rightarrow 2^{X},$ $Q:T(X\times X)\times X\times X\rightarrow
2^{Z}$\textit{\ and }$G:T(X\times X)\times X\times X\rightarrow 2^{Z}$ 
\textit{be correspondences with nonempty values. Assume that:}
\end{theorem}

\textit{i) }$S_{1}^{-1}$\textit{\ and }$S_{2}^{-1}$\textit{\ are open-valued
and for each }$u\in X,$\textit{\ the set }$\{x\in X:\exists t\in
T(x,u),F(t,u,x)\cap G(t,x,x)=\emptyset \}$\textit{\ is open; }

\textit{ii)} $S_{1}$\textit{\ and }$S_{2}$\textit{\ are convex valued and
for each }$x\in X,$\textit{\ the set }$\{u\in X:\exists t\in
T(x,u),F(t,u,x)\cap G(t,x,x)=\emptyset \}$\textit{\ is convex; }

\textit{iii)} \textit{the set }$A=$ $\{x\in X:$\textit{\ there exist }$u\in
S_{2}(x)$ \textit{and }$t\in T(x,u)$\textit{\ such that} $F(t,u,x)\cap
G(t,x,x)=\emptyset \}$ is closed;

\textit{iv) }$\mathit{\{x}\in X:x\in S_{2}(x)\}=\emptyset ;$

\textit{v) there exists }$M\in \langle A\rangle $\textit{\ such that }$%
\bigcup\nolimits_{u\in M}[\{x\in X:\exists t\in T(x,u),F(t,u,x)\cap
G(t,x,x)=\emptyset \}=X$ or $\bigcup\nolimits_{u\in M}[(X\backslash A)\cap 
\mathit{\ }S_{1}^{-1}(u)]=X.$

\textit{Then, there exists }$x^{\ast }\in S_{1}(x^{\ast })$\textit{\ such
that }$\forall y\in S_{2}(x^{\ast }),$ $\forall t^{\ast }\in T(x^{\ast },y),$
$F(t^{\ast },y,x^{\ast })\cap G(t^{\ast },x^{\ast },x^{\ast })\neq \emptyset 
$\textit{.\medskip }

\textit{Proof.} Let $P,G:X\rightarrow 2^{X}$ be defined by

$P(x)=\{u\in X:\exists t\in T(x,u),F(t,u,x)\cap G(t,x,x)=\emptyset \},$ for
each $x\in X$ and

$G(x)=S_{2}(x)\cap P(x),$ for each $x\in X.$

We are going to show that there exists $x^{\ast }\in X$ such that $x^{\ast
}\in S_{1}(x^{\ast })$ and $S_{2}(x^{\ast })\cap P(x^{\ast })=\emptyset .$

The rest of the proof follows the same line as the proof of Theorem
16.\medskip

We obtain a new existence theorem of solutions for a generalized vector
variational inclusion problem.

\begin{theorem}
\textit{Let }$Z$\textit{\ be a Hausdorff topological vector space and let }$%
X $\textit{\ be a nonempty compact convex subset of a\ topological vector
space }$E$\textit{. Let }$S_{1},S_{2}:X\rightarrow 2^{X},$\textit{\ }$%
T:X\times X\rightarrow 2^{X},$ $F:T(X\times X)\times X\times X\rightarrow
2^{Z}$\textit{\ and }$G:T(X\times X)\times X\times X\rightarrow 2^{Z}$ 
\textit{be correspondences with nonempty values. Assume that:}
\end{theorem}

\textit{i) }$S_{2}^{-1}$\textit{\ is open-valued and for each }$u\in X,$%
\textit{\ the set }$\{x\in X:\exists t\in T(x,u),F(t,u,x)\nsubseteq
G(t,x,x)\}$\textit{\ is open; }

\textit{ii)} $S_{2}$\textit{\ is convex valued and for each }$x\in X,$%
\textit{\ the set }$\{u\in X:\exists t\in T(x,u),F(t,u,x)\nsubseteq
G(t,x,x)\}$\textit{\ is convex; }

\textit{iii) if }$A=\{x\in X:\mathit{\ }$\textit{there exist }$u\in S_{2}(x)$
\textit{and} $t\in T(x,u)$ \textit{such that} $F(t,u,x)\nsubseteq
G(t,x,x)\}, $\textit{\ then,} \textit{for each }$N\in \langle A\rangle ,$%
\textit{\ }$($\textit{co}$N\smallsetminus N)\cap A=\emptyset ;$

\textit{iv) }$\mathit{\{x}\in X:x\in S_{2}(x)\}=\emptyset ;$

\textit{v) there exists }$M\in \langle A\rangle $\textit{\ such that }$%
\bigcup\nolimits_{u\in M}[\{x\in X:\exists t\in T(x,u),F(t,u,x)\nsubseteq
G(t,x,x)\}\cap S_{2}^{-1}(u)]=X$ or $\bigcup\nolimits_{u\in M}[(X\backslash
A)\cap \mathit{\ }S_{1}^{-1}(u)]=X.$

\textit{Then, there exists }$x^{\ast }\in S_{1}(x^{\ast })$\ \textit{such
that }$\forall y\in S_{2}(x^{\ast }),$ $\forall t^{\ast }\in T(x^{\ast },y),$
$F(t^{\ast },y,x^{\ast })\subseteq G(t^{\ast },x^{\ast },x^{\ast })$\textit{%
.\medskip }

\textit{Proof.} Let $P,G:X\rightarrow 2^{X}$ be defined by

$P(x)=\{u\in X:\exists t\in T(x,u),F(t,u,x)\nsubseteq G(t,x,x)\},$ for each $%
x\in X$ and

$G(x)=S_{2}(x)\cap P(x),$ for each $x\in X.$

We are going to show that there exists $x^{\ast }\in X$ such that $x^{\ast
}\in S_{1}(x^{\ast })$ and $S_{2}(x^{\ast })\cap P(x^{\ast })=\emptyset .$

We consider two cases.

Case I.

$A=\{x\in X:$\textit{\ }$P(x)\cap K(x)\neq \emptyset \}=\{x\in X:$\textit{\ }%
$G(x)\neq \emptyset \}$\textit{\ }is nonempty.

The correspondence $G:A\rightarrow 2^{X},$ defined by $G(x)=S_{2}(x)\cap
P(x) $ for each $x\in A,$ is nonempty valued on $A.$ We note that for each $%
u\in X,$ $G^{-1}(u)=P^{-1}(u)\cap S_{2}^{-1}(u)$ is a convex set since it is
an intersection of convex sets.

Further, let us define the correspondences $H,L:X\rightarrow 2^{X}$ by

$H(x)=\left\{ 
\begin{array}{c}
G(x),\text{ if }x\in A; \\ 
\emptyset ,\text{\ \ \ \ otherwise}%
\end{array}%
\right. =\left\{ 
\begin{array}{c}
S_{2}(x)\cap P(x),\text{ if }x\in A; \\ 
\emptyset ,\text{ \ \ \ \ \ \ \ \ \ \ \ \ \ \ \ \ otherwise.}%
\end{array}%
\right. $

According to i) and ii), $H$ is open convex valued.

$L(x)=\left\{ 
\begin{array}{c}
G(x),\text{ \ \ \ \ \ \ \ \ \ \ \ if }x\in A; \\ 
S_{1}(x),\text{ \ \ \ \ \ \ \ \ \ otherwise.}%
\end{array}%
\right. $

For each $u\in X,$

$H^{-1}(u)=\{x\in X:$\textit{\ }$\mathit{u}\in H(x)\}=$

\ \ \ \ \ \ \ \ \ \ \ $=\{x\in A:$\textit{\ }$\mathit{u}\in G(x)\}=$

\ \ \ \ \ \ \ \ \ \ \ $=\{x\in A:$\textit{\ }$\mathit{u}\in S_{2}(x)\cap
P(x)\}$

$\ \ \ \ \ \ \ \ =A\cap P^{-1}(u)\cap S_{2}^{-1}(u)=$

\ \ \ \ \ \ \ \ \ $=P^{-1}(u)\cap S_{2}^{-1}(u).$

Since for each $u\in X,$ $S_{2}^{-1}(u)$ and $P^{-1}(u)$ are open, then, $%
H^{-1}(u)$ is open.

Assumption v) implies that there exists $M\in \langle A\rangle $ such that $%
\tbigcup\nolimits_{x\in M}H^{-1}(x)=X.$

Let us define $Q:X\rightarrow 2^{X}$ by $Q(x):=X\backslash H^{-1}(x)$ for
each\textit{\ }$x\in X.$

Then, $Q$ is closed-valued and $\tbigcap\nolimits_{x\in M}Q(x)=X\backslash
\tbigcup\nolimits_{x\in M}H^{-1}(x)=\emptyset .$

According to Lemma 4, we can conclude that $Q$ is not a KKM correspondence.
Thus, there exists $N\in \langle X\rangle $ such that\textit{\ }co$%
N\varsubsetneq Q(N)=\tbigcup\nolimits_{x\in N}(X\backslash H^{-1}(x)).$

Hence, there exists $x^{\ast }\in $co$N$ with the property that $x^{\ast
}\in H^{-1}(x)$ for each $x\in N.$ Therefore, there exists $x^{\ast }\in $co$%
N$ such that $x^{\ast }\in H^{-1}(x)$ for each $x\in N,$ which implies $%
N\subset H(x^{\ast })$. Further, it is true that co$N\subset $co$H(x^{\ast
})\subset L(x^{\ast })$. Consequently, $x^{\ast }\in $co$N\subset $co$%
H(x^{\ast })\subset L(x^{\ast }),$ which means that $x^{\ast }\in L(x^{\ast
}),$ that is, $x^{\ast }\in $co$N$ is a fixed point for $L.$

We notice that, according to $iv)$, $x\notin S_{2}(x)$ for each $x\in X,$
and then, $x^{\ast }\notin A.$ This fact is possible since $x^{\ast }\in $co$%
N$ and assumption iii) asserts that $($co$N\smallsetminus N)\cap A=\emptyset
.$ Therefore, $x^{\ast }\in S_{1}(x^{\ast })$ and $G(x^{\ast
})=S_{2}(x^{\ast })\cap P(x^{\ast })=\emptyset .$

Consequently, there exists $x^{\ast }\in X$ such that $x^{\ast }\in
S_{1}(x^{\ast })$ and $\forall y\in S_{2}(x^{\ast }),$ $\forall t^{\ast }\in
T(x^{\ast },y),$ $F(t^{\ast },y,x^{\ast })\subseteq G(t^{\ast },x^{\ast
},x^{\ast })$\textit{.}

Case II.

$A=\{x\in X:$\textit{\ }$P(x)\cap S_{2}(x)\neq \emptyset \}=\{x\in X:$%
\textit{\ }$G(x)\neq \emptyset \}\mathit{=}\emptyset $.

In this case, $G(x)=\emptyset $ for each $x\in X.$

Let us define $Q:X\rightarrow 2^{X}$ by $Q(x):=X\backslash S_{1}(x)$ for each%
\textit{\ }$x\in X.$ The proof follows the same line as above and we obtain
that there exists $x^{\ast }\in X$ such that $x^{\ast }\in S_{1}(x^{\ast }).$

Obviously, $G(x^{\ast })=\emptyset .$ Consequently, the conclusion holds in
Case II.

\begin{remark}
Assumption i) can be replaced with
\end{remark}

i') $S_{2}^{-1}$\textit{\ }is closed-valued and for each $u\in X,$\ the set $%
\{x\in X:\exists t\in T(x,u),F(t,u,x)\nsubseteq G(t,x,x)\}$\ is closed.\ 

In this case, $Q$ is open-valued.\smallskip

Theorem 19 can be stated as follows.

\begin{theorem}
\textit{Let }$Z$\textit{\ be a Hausdorff topological vector space and let }$%
X $\textit{\ be a nonempty compact convex subset of a\ topological vector
space }$E$\textit{. Let }$S_{1},S_{2}:X\rightarrow 2^{X},$\textit{\ }$%
T:X\times X\rightarrow 2^{X},$ $F:T(X\times X)\times X\times X\rightarrow
2^{Z}$\textit{\ and }$G:T(X\times X)\times X\times X\rightarrow 2^{Z}$ 
\textit{be correspondences with nonempty values. Assume that:}
\end{theorem}

\textit{i) }$S_{2}^{-1}$\textit{\ is open-valued and for each }$u\in X,$%
\textit{\ the set }$\{x\in X:\exists t\in T(x,u),F(t,u,x)\cap
G(t,x,x)=\emptyset \}$\textit{\ is open; }

\textit{ii)} $S_{2}$\textit{\ is convex valued and for each }$x\in X,$%
\textit{\ the set }$\{u\in X:\exists t\in T(x,u),F(t,u,x)\cap
G(t,x,x)=\emptyset \}$\textit{\ is convex; }

\textit{iii) if }$A=$ $\{x\in X:$\textit{\ there exist }$u\in S_{2}(x)$ 
\textit{and }$t\in T(x,u)$\textit{\ such that} $F(t,u,x)\cap
G(t,x,x)=\emptyset \},$ \textit{then, for each }$N\in \langle A\rangle ,$%
\textit{\ }$($\textit{co}$N\smallsetminus N)\cap A=\emptyset ;$

\textit{iv) }$\mathit{\{x}\in X:x\in S_{2}(x)\}=\emptyset ;$

\textit{v) there exists }$M\in \langle A\rangle $\textit{\ such that }$%
\bigcup\nolimits_{u\in M}[\{x\in X:\exists t\in T(x,u),F(t,u,x)\cap
G(t,x,x)=\emptyset \}\cap S_{2}^{-1}(u)]=X$ or $\bigcup\nolimits_{u\in
M}[(X\backslash A)\cap \mathit{\ }S_{1}^{-1}(u)]=X.$

\textit{Then, there exists }$x^{\ast }\in S_{1}(x^{\ast })$\ \textit{such
that }$\forall y\in S_{2}(x^{\ast }),$ $\forall t^{\ast }\in T(x^{\ast },y),$
$F(t^{\ast },y,x^{\ast })\cap G(t^{\ast },x^{\ast },x^{\ast })=\emptyset $%
\textit{.\medskip }

\textit{Proof.} Let $P,G:X\rightarrow 2^{X}$ be defined by

$P(x)=\{u\in X:\exists t\in T(x,u),F(t,u,x)\cap G(t,x,x)=\emptyset \},$ for
each $x\in X$ and

$G(x)=S_{2}(x)\cap P(x),$ for each $x\in X.$

We are going to show that there exists $x^{\ast }\in X$ such that $x^{\ast
}\in S_{1}(x^{\ast })$ and $S_{2}(x^{\ast })\cap P(x^{\ast })=\emptyset .$

The rest of the proof follows a similar line as the proof of Theorem
18.\medskip

\begin{theorem}
\textit{Let }$Z$\textit{\ be a Hausdorff topological vector space and let }$%
X $\textit{\ be a nonempty compact convex subset of a\ topological vector
space }$E$\textit{. Let }$S_{1},S_{2}:X\rightarrow 2^{X},$\textit{\ }$%
T:X\times X\rightarrow 2^{X}$,\textit{\ }$F:T(X\times X)\times X\times
X\rightarrow 2^{Z}$\textit{and }$G:T(X\times X)\times X\times X\rightarrow
2^{Z}$\textit{\ be correspondences with nonempty values. Assume that:}
\end{theorem}

\textit{i) }$S_{2}$\textit{\ is open-valued and for each }$x\in X,$\textit{\
the set }$\{u\in X:\exists t\in T(x,u),F(t,u,x)\nsubseteq G(t,x,x)\}$\textit{%
\ is open; }

\textit{ii) if }$A=\{x\in X:\mathit{\ }$\textit{there exist }$\mathit{u}\in 
\mathit{S}_{2}\mathit{(x)}$\textit{\ and }$\mathit{t}\in \mathit{T(x,u)}$%
\textit{\ such that }$F(t,u,x)\nsubseteq G(t,x,x)\},$ \textit{then,} \textit{%
for each }$N\in \langle A\rangle ,$\textit{\ }$($\textit{co}$N\smallsetminus
N)\cap A=\emptyset ;$

\textit{iii)} $\mathit{\{x}\in X:x\in S_{2}(x)\}=\emptyset $\textit{;}

\textit{iv) there exists }$M\in \langle A\rangle $\textit{\ such that }$%
\bigcup\nolimits_{x\in M}[\{u\in X:\exists t\in T(x,u),F(t,u,x)\nsubseteq
G(t,x,x)\}\cap S_{2}(x)]=X;$

\textit{v) }$S_{2}^{-1}:X\rightarrow 2^{X}$\textit{\ is convex valued and
for each }$u\in X,$\textit{\ the set }$\{x\in X:\exists t\in
T(x,u),F(t,u,x)\nsubseteq G(t,x,x)\}$\textit{\ is convex;}

\textit{Then, there exists }$x^{\ast }\in S_{1}(x^{\ast })$\ \textit{such
that }$\forall y\in S_{2}(x^{\ast }),$ $\forall t^{\ast }\in T(x^{\ast },y),$
$F(t^{\ast },y,x^{\ast })\subseteq G(t^{\ast },x^{\ast },x^{\ast })$\textit{%
.\medskip }

\textit{Proof.} Let $P,G:X\rightarrow 2^{X}$ be defined by

$P(x)=\{u\in X:\exists t\in T(x,u),F(t,u,x)\nsubseteq G(t,x,x)\},$ for each $%
x\in X$ and

$G(x)=S_{2}(x)\cap P(x),$ for each $x\in X.$

We are going to show that there exists $x^{\ast }\in X$ such that $x^{\ast
}\in S_{1}(x^{\ast })$ and $S_{2}(x^{\ast })\cap P(x^{\ast })=\emptyset .$

We consider two cases.

Case I.

$A=\{x\in X:$\textit{\ }$S_{2}(x)\cap P(x)\neq \emptyset \}=\{x\in X:$%
\textit{\ }$G(x)\neq \emptyset \}$\textit{\ }is nonempty.

The correspondence $G:A\rightarrow 2^{X},$ defined by $G(x)=S_{2}(x)\cap
P(x) $ for each $x\in A,$ is nonempty valued on $A.$ We note that for each $%
u\in X,$ $G^{-1}(u)=P^{-1}(u)\cap S_{2}^{-1}(u)$ is a convex set as
intersection of convex sets.

Further, let us define the correspondences $H,L:X\rightarrow 2^{X}$ by

$H(x)=\left\{ 
\begin{array}{c}
G(x),\text{ if }x\in A; \\ 
\emptyset ,\text{\ \ \ \ otherwise}%
\end{array}%
\right. $

$L(x)=\left\{ 
\begin{array}{c}
G(x),\text{ if }x\in A; \\ 
S_{1}(x),\text{\ otherwise}%
\end{array}%
\right. =\left\{ 
\begin{array}{c}
S_{2}(x)\cap P(x),\text{ if }x\in A; \\ 
S_{1}(x),\text{ \ \ \ \ \ \ \ \ \ otherwise.}%
\end{array}%
\right. $

According to i), $H$ is open-valued and according to v),for each $u\in X,$ $%
P^{-1}(u)$ is convex.

For each $u\in X,$

$H^{-1}(u)=\{x\in X:$\textit{\ }$\mathit{u}\in H(x)\}=$

\ \ \ \ \ \ \ \ \ \ \ $=\{x\in A:$\textit{\ }$\mathit{u}\in G(x)\}=$

\ \ \ \ \ \ \ \ \ \ \ $=\{x\in A:$\textit{\ }$\mathit{u}\in S_{2}(x)\cap
P(x)\}$

$\ \ \ \ \ \ \ \ =A\cap P^{-1}(u)\cap S_{2}^{-1}(u)=$

\ \ \ \ \ \ \ \ \ $=P^{-1}(u)\cap S_{2}^{-1}(u).$

$L^{-1}(u)=\{x\in X:$\textit{\ }$\mathit{u}\in L(x)\}=$

\ \ \ \ \ \ \ \ \ \ \ $=\{x\in A:$\textit{\ }$\mathit{u}\in G(x)\}\cup
\{x\in X\setminus A:$\textit{\ }$\mathit{u}\in S_{1}(x)\}=$

\ \ \ \ \ \ \ \ \ \ \ $=\{x\in A:$\textit{\ }$\mathit{u}\in S_{2}(x)\cap
P(x)\}\cup \{x\in X\setminus A:\mathit{\ u}\in S_{1}(x)\}$

$\ \ \ \ \ \ \ \ =(A\cap P^{-1}(u)\cap S_{2}^{-1}(u)\cup \lbrack
(X\backslash A)\cap $\textit{\ }$S_{1}^{-1}(u)]=$

\ \ \ \ \ \ \ $=[P^{-1}(u)\cap S_{2}^{-1}(u)]\cup \lbrack (X\backslash
A)\cap $\textit{\ }$S_{1}^{-1}(u)]$

Since for each $u\in X,$ $S_{2}^{-1}(u)$ and $P^{-1}(u)$ are convex, then, $%
H^{-1}(u)$ is convex. Therefore, co$H^{-1}(u)\subset L^{-1}(u)$ for each $%
u\in X$.

Assumption iv) implies that there exists $M\in \langle A\rangle $ such that $%
\tbigcup\nolimits_{x\in M}H(x)=X.$

Let us define $Q:X\rightarrow 2^{X}$ by $Q(x):=X\backslash H(x)$ for each%
\textit{\ }$x\in X.$

Then, $Q$ is closed-valued and $\tbigcap\nolimits_{x\in M}Q(x)=X\backslash
\tbigcup\nolimits_{x\in M}H(x)=\emptyset .$

According to Lemma 4, we can conclude that $Q$ is not a KKM correspondence.
Thus, there exists $N\in \langle X\rangle $ such that\textit{\ }co$%
N\varsubsetneq Q(N)=\tbigcup\nolimits_{x\in N}(X\backslash H(x)).$

Hence, there exists $x^{\ast }\in $co$N$ with the property that $x^{\ast
}\in H(x)$ for each $x\in N.$ Therefore, there exists $x^{\ast }\in $co$N$
such that $x^{\ast }\in H(x)$ for each $x\in N,$ which implies $N\subset
H^{-1}(x^{\ast })$. Further, it is true that co$N\subset $co$H^{-1}(x^{\ast
})\subset L^{-1}(x^{\ast })$. Consequently, $x^{\ast }\in $co$N\subset $co$%
H^{-1}(x^{\ast })\subset L^{-1}(x^{\ast }),$ which means that $x^{\ast }\in
L(x^{\ast }),$ that is, $x^{\ast }\in $co$N$ is a fixed point for $L.$

We notice that, according to iii), $x\notin S_{2}(x)$ for each $x\in X,$ and
then, $x^{\ast }\notin A.$ This fact is possible since $x^{\ast }\in $co$N$
and assumption iii) asserts that $($co$N\smallsetminus N)\cap A=\emptyset .$
Therefore, $x^{\ast }\in S_{1}(x^{\ast })$ and $G(x^{\ast })=S_{2}(x^{\ast
})\cap P(x^{\ast })=\emptyset .$

Consequently, there exists $x^{\ast }\in S_{1}(x^{\ast })$\ such that\textit{%
\ }$\forall y\in S_{2}(x^{\ast }),$ $\forall t^{\ast }\in T(x^{\ast },y),$ $%
F(t^{\ast },y,x^{\ast })\subseteq G(t^{\ast },x^{\ast },x^{\ast })$\textit{.}

Case II.

$A=\{x\in X:$\textit{\ }$P(x)\cap S_{2}(x)\neq \emptyset \}=\{x\in X:$%
\textit{\ }$G(x)\neq \emptyset \}\mathit{=}\emptyset $.

In this case, $G(x)=\emptyset $ for each $x\in X.$

Let us define $Q:X\rightarrow 2^{X}$ by $Q(x):=X\backslash S_{1}(x)$ for each%
\textit{\ }$x\in X.$ The proof follows the same line as above and we obtain
that there exists $x^{\ast }\in X$ such that $x^{\ast }\in S_{1}(x^{\ast }).$

Obviously, $G(x^{\ast })=\emptyset .$ Consequently, the conclusion holds in
Case II.

\begin{remark}
Assumption i) can be replaced with
\end{remark}

i') $S_{2}$\ is closed-valued and for each $x\in X,$ the set $\{u\in
X:\exists t\in T(x,u),F(t,u,x)\nsubseteq G(t,x,x)\}$ is closed.\textit{\ }

In this case, $Q$ is open-valued.\smallskip

Now, we are proving the existence of solutions for a general vector
variational inclusion problem concerning correspondences under new
assumptions.

\begin{theorem}
\textit{Let }$Z$\textit{\ be a Hausdorff topological vector space and let }$%
X $\textit{\ be a nonempty compact convex subset of a\ topological vector
space }$E$\textit{. Let }$S_{1},S_{2}:X\rightarrow 2^{X},$\textit{\ }$%
T:X\times X\rightarrow 2^{X}$ \textit{and }$F:T(X\times X)\times X\times
X\rightarrow 2^{Z}$\textit{\ and }$G:T(X\times X)\times X\times X\rightarrow
2^{Z}$ \textit{be correspondences with nonempty values. Assume that:}
\end{theorem}

\textit{i) }$S_{1},S_{2}$\textit{\ are open-valued and for each }$(x,y)\in
X\times X$\textit{\ with the property that }$\exists t\in T(x,y)$ such that $%
F(t,y,x)\nsubseteq G(t,x,x)\},$\textit{\ there exists }$z=z_{x,y}\in X$%
\textit{\ such that }$y\in $\textit{int}$_{X}\{u\in X:\exists t\in
T(z_{x,y},u),F(t,u,z_{x,y})\nsubseteq G(t,z_{x,y},z_{x,y})\}$\textit{; }

\textit{ii) there exists }$M\in \langle X\rangle $\textit{\ such that}%
\newline
$\bigcup\nolimits_{x\in M\cap A}[\bigcup\nolimits_{\{y\in S_{2}(x),\text{ }%
t\in T(x,y):\text{ }F(t,y,x)\nsubseteq G(t,x,x)\}}($\textit{int}$_{X}\{u\in
X:\exists t\in T(z_{x,y},u),$ \newline
$F(t,u,z_{x,y})\nsubseteq G(t,z_{x,y},z_{x,y})\}\cap S_{2}(x)]\bigcup
\bigcup\nolimits_{x\in M\backslash A}S_{1}(x)=X,$\textit{\ where }\newline
$A=$ $\{x\in X:$\textit{\ there exist }$u\in S_{2}(x)$,\textit{\ }$t\in
T(x,u)$\textit{\ such that} $F(t,u,x)\nsubseteq G(t,x,x)\}$\textit{;}

\textit{iii)} $\mathit{\{x}\in X:x\in S_{2}(x)\}=\emptyset $\textit{;}

\textit{iv) for each }$u\in X,$\textit{\ the set }$M_{u}$\textit{\ is
convex, where}\newline
$M_{u}=\{x\in X:\exists t\in T(x,u),F(t,u,x)\nsubseteq G(t,x,x)\}\cup
\lbrack (X\backslash A)\cap (S_{1}^{-1}(u)];$

\textit{v) for each }$x\in A,$\textit{\ }\newline
$\bigcup\nolimits_{\{y\in S_{2}(x),\text{ }t\in T(x,y):\text{ }%
F(t,y,x)\nsubseteq G(t,x,x)\}}($\textit{int}$_{X}\{u\in X:\exists t\in
T(z_{x,y},u),F(t,u,z_{x,y})\nsubseteq G(t,z_{x,y},z_{x,y})\}\subset \{u\in
X:\exists t\in T(x,u),F(t,u,x)\nsubseteq G(t,x,x)\}$\textit{;}

\textit{Then, there exists }$x^{\ast }\in S_{1}(x^{\ast })$\ \textit{such
that }$\forall y\in S_{2}(x^{\ast }),$ $\forall t^{\ast }\in T(x^{\ast },y),$
$F(t^{\ast },y,x^{\ast })\subseteq G(t^{\ast },x^{\ast },x^{\ast })$.\textit{%
\medskip }

\textit{Proof.} Let $P,G:X\rightarrow 2^{X}$ be defined by

$P(x)=\{u\in X:\exists t\in T(x,u),F(t,u,x)\nsubseteq G(t,x,x)\},$ for each $%
x\in X$ and

$G(x)=S_{2}(x)\cap P(x),$ for each $x\in X.$

We are going to show that there exists $x^{\ast }\in X$ such that $x^{\ast
}\in S_{1}(x^{\ast })$ and $S_{2}(x^{\ast })\cap P(x^{\ast })=\emptyset .$

We consider two cases.

Case I.

$A=\{x\in X:$\textit{\ }$P(x)\cap S_{2}(x)\neq \emptyset \}=\{x\in X:$%
\textit{\ }$G(x)\neq \emptyset \}$\textit{\ }is nonempty.

The correspondence $G:A\rightarrow 2^{X},$ defined by $G(x)=S_{2}(x)\cap
P(x) $ for each $x\in A,$ is nonempty on $A.$

Further, let us define the correspondences $H,L:X\rightarrow 2^{X}$ by

$H(x)=\left\{ 
\begin{array}{c}
G(x),\text{ if }x\in A; \\ 
S_{1}(x),\text{\ otherwise}%
\end{array}%
\right. $ and

$L(x)=\left\{ 
\begin{array}{c}
P(x),\text{ if }x\in A; \\ 
S_{1}(x),\text{\ otherwise}%
\end{array}%
\right. $

For each $u\in X,$

$L^{-1}(u)=\{x\in X:$\textit{\ }$\mathit{u}\in L(x)\}=$

\ \ \ \ \ \ \ \ \ \ \ $=\{x\in A:$\textit{\ }$\mathit{u}\in P(x)\}\cup
\{x\in X\setminus A:$\textit{\ }$\mathit{u}\in S_{1}(x)\}=$

$\ \ \ \ \ \ \ \ =(A\cap P^{-1}(u))\cup \lbrack (X\backslash A)\cap $\textit{%
\ }$S_{1}^{-1}(u)]=$

\ \ \ \ \ \ \ \ \ $=P^{-1}(u)\cup \lbrack (X\backslash A)\cap
S_{1}^{-1}(u)]. $

According to i), $G$ is transfer open-valued. Assumptions i) and ii) imply
that there exists $M\in \langle X\rangle $ and for each $x\in M$ and $y\in
H(x),$ there exists $z_{x,y}\in X$ such that $y\in $int$_{X}H(z_{x,y})\cap
H(x)$ and $\tbigcup\nolimits_{x\in M}(\tbigcup\nolimits_{y\in H(x)}$int$%
_{X}H(z_{x,y}))=X$. In addition, assumption v) implies $\tbigcup\nolimits_{y%
\in H(x)}$int$_{X}H(z_{x,y})\subseteq P(x)$ for each $x\in A.$ We note that
if $x\in X\backslash A,$ then $H(x)=S_{1}(x)$ is open and $y\in H(x)$
implies $z_{x,y\text{ }}=x$ and $y\in H(z_{x.y})=$int$H(x).$ In this last
case it is obvious that $\tbigcup\nolimits_{y\in H(x)}$int$%
_{X}H(z_{x,y})=\tbigcup\nolimits_{y\in H(x)}$int$_{X}H(x)=\tbigcup%
\nolimits_{y\in H(x)}H(x)=H(x).$

Let us define $Q:X\rightarrow 2^{X}$ by $Q(x):=X\backslash
\tbigcup\nolimits_{y\in H(x)}$int$_{X}H(z_{x,y})$ for each\textit{\ }$x\in
H. $

Then, $Q$ is closed-valued and $\tbigcap\nolimits_{x\in M}Q(x)=X\backslash
\tbigcup\nolimits_{x\in M}(\tbigcup\nolimits_{y\in H(x)}$int$%
_{X}H(z_{x,y}))=\emptyset .$

According to Lemma 4, we can conclude that $Q$ is not a KKM correspondence.
Thus, there exists $N\in \langle X\rangle $ such that\textit{\ }co$%
N\varsubsetneq Q(N)=\tbigcup\nolimits_{x\in N}(X\backslash
\tbigcup\nolimits_{y\in H(x)}$int$_{X}H(z_{x,y})).$

Hence, there exists $x^{\ast }\in $co$N$ with the property that $x^{\ast
}\in \tbigcup\nolimits_{y\in H(x)}$int$_{X}H(z_{x,y})$ for each $x\in N.$ If 
$x\in A,$ $\tbigcup\nolimits_{y\in H(x)}$int$_{X}H(z_{x,y})\subset P(x)$ and
if $x\in X\backslash A,$ $\tbigcup\nolimits_{y\in H(x)}$int$%
_{X}H(z_{x,y})=H(x).$ Therefore, there exists $x^{\ast }\in $co$N$ such that 
$x^{\ast }\in L(x)$ for each $x\in N,$ which implies $N\subset
L^{-1}(x^{\ast })$. Further, it is true that co$N\subset $co$L^{-1}(x^{\ast
})=L^{-1}(x^{\ast })$. Consequently, $x^{\ast }\in $co$N\subset $co$%
L^{-1}(x^{\ast })=L^{-1}(x^{\ast }),$ which means that $x^{\ast }\in
L(x^{\ast }).$ We notice that, according to ii), $x\notin S_{2}(x)$ for each 
$x\in X,$ and then, $x^{\ast }\notin A.$ Therefore, $x^{\ast }\in
S_{1}(x^{\ast })$ and $G(x^{\ast })=S_{2}(x^{\ast })\cap P(x^{\ast
})=\emptyset .$

Consequently, there exists $x^{\ast }\in S_{1}(x^{\ast })$\ such that $%
\forall y\in S_{2}(x^{\ast }),$ $\forall t^{\ast }\in T(x^{\ast },y),$ $%
F(t^{\ast },y,x^{\ast })\subseteq G(t^{\ast },x^{\ast },x^{\ast }).$

Case II.

$A=\{x\in X:$\textit{\ }$P(x)\cap S_{2}(x)\neq \emptyset \}=\{x\in X:$%
\textit{\ }$G(x)\neq \emptyset \}\mathit{=}\emptyset $.

In this case, $G(x)=\emptyset $ for each $x\in X.$ Let us define $%
Q:X\rightarrow 2^{X}$ by $Q(x):=X\backslash S_{1}(x)$ for each\textit{\ }$%
x\in X.$ The proof follows the same line as above and we obtain that there
exists $x^{\ast }\in X$ such that $x^{\ast }\in S_{1}(x^{\ast }).$

Obviously, $G(x^{\ast })=\emptyset .$ Consequently, the conclusion holds in
Case II.\smallskip

Now, we are establishing Theorem 22.

\begin{theorem}
\textit{Let }$Z$\textit{\ be a Hausdorff topological vector space and let }$%
X $\textit{\ be a nonempty compact convex subset of a\ topological vector
space }$E$\textit{. Let }$S_{1},S_{2}:X\rightarrow 2^{X},$\textit{\ }$%
T:X\times X\rightarrow 2^{X},$ $F:T(X\times X)\times X\times X\rightarrow
2^{Z}$\textit{\ and }$G:T(X\times X)\times X\times X\rightarrow 2^{Z}$ 
\textit{be correspondences with nonempty values. Assume that:}
\end{theorem}

\textit{i) }$S_{1},S_{2}$\textit{\ are open-valued and for each }$(x,y)\in
X\times X$\textit{\ with the property that }$\exists t\in T(x,y)$ such that $%
F(t,y,x)\cap G(t,x,x)=\emptyset \},$\textit{\ there exists }$z=z_{x,y}\in X$%
\textit{\ such that }$y\in $\textit{int}$_{X}\{u\in X:\exists t\in
T(z_{x,y},u),F(t,u,z_{x,y})\cap G(t,z_{x,y},z_{x,y})=\emptyset \}$\textit{; }

\textit{ii) there exists }$M\in \langle X\rangle $\textit{\ such that}

$\bigcup\nolimits_{x\in M\cap A}[\bigcup\nolimits_{\{y\in S_{2}(x),\text{ }%
t\in T(x,y)\text{: }F(t,y,x)\nsubseteq G(t,x,x)\}}($\textit{int}$_{X}\{u\in
X:\exists t\in T(z_{x,y},u),$ \newline
$F(t,u,z_{x,y})\cap G(t,z_{x,y},z_{x,y})=\emptyset \}\cap S_{2}(x)]\bigcup
\bigcup\nolimits_{x\in M\backslash A}S_{1}(x)=X,$\textit{\ where }\newline
$A=$ $\{x\in X:$\textit{\ there exist }$u\in S_{2}(x)$,\textit{\ }$t\in
T(x,u)$\textit{\ such that} $F(t,u,x)\cap G(t,x,x)=\emptyset \}$\textit{;}

\textit{iii) }$\mathit{\{x}\in X:x\in S_{2}(x)\}=\emptyset $\textit{; }

\textit{iv) for each }$u\in X,$\textit{\ the set }$M_{u}$\textit{\ is
convex, where}\newline

$M_{u}=\{x\in X:\exists t\in T(x,u),F(t,u,x)\nsubseteq G(t,x,x)\}\cup
\lbrack (X\backslash A)\cap (S_{1}^{-1}(u)];$

\textit{v) for each }$x\in A,$\textit{\ }$\bigcup\nolimits_{y\in S_{2}(x),%
\text{ }\exists t\in T(x,y)\text{ such that }F(t,y,x)\nsubseteq G(t,x,x)}($%
\textit{int}$_{X}\{u\in X:\exists t\in T(z_{x,y},u),F(t,u,z_{x,y})\cap
G(t,z_{x,y},z_{x,y})=\emptyset \}\subset $

$\{u\in X:\exists t\in T(x,u),F(t,u,x)\cap G(t,x,x)=\emptyset \}$\textit{;}

\textit{Then, there exists }$x^{\ast }\in S_{1}(x^{\ast })$\ \textit{such
that }$\forall y\in S_{2}(x^{\ast }),$ $\forall t^{\ast }\in T(x^{\ast },y),$
$F(t^{\ast },y,x^{\ast })\cap G(t^{\ast },x^{\ast },x^{\ast })\neq \emptyset 
$.\textit{\medskip }

We establish sufficient conditions which assure the existence of solutions
for a general vector variational inclusion problem.

\begin{theorem}
\textit{Let }$Z$\textit{\ be a Hausdorff topological vector space and let }$%
X $\textit{\ be a nonempty compact convex subset of a\ topological vector
space }$E$\textit{. Let }$S_{1},S_{2}:X\rightarrow 2^{X},$\textit{\ }$%
T:X\times X\rightarrow 2^{X}$, $F:T(X\times X)\times X\times X\rightarrow
2^{Z}$\textit{\ and }$G:T(X\times X)\times X\times X\rightarrow 2^{Z}$ 
\textit{be correspondences with nonempty values. Assume that:}
\end{theorem}

\textit{i) }$S_{1},S_{2}$\textit{\ are open-valued and for each }$(x,y)\in
X\times X$\textit{\ with the property that }$\exists t\in T(x,y)$ such that $%
F(t,y,x)\nsubseteq G(t,x,x)\},$\textit{\ there exists }$z=z_{x,y}\in X$%
\textit{\ such that }$y\in $\textit{int}$_{X}\{u\in X:\exists t\in
T(z_{x,y},u),F(t,u,z_{x,y})\nsubseteq G(t,z_{x,y},z_{x,y})\}$\textit{; }

\textit{ii) for each }$x\in X$\textit{\ and} $t\in T(x,x),F(t,x,x)\subseteq
G(t,x,x);$\textit{\ }

\textit{iii)} \textit{if} $A=$ $\{x\in X:$\ \textit{there exist }$u\in
S_{2}(x)$ \textit{and }$t\in T(x,u)$\textit{\ such that} $F(t,u,x)\nsubseteq
G(t,x,x)\},$\textit{\ then for each }$N\in \langle A\rangle ,$\textit{\ }$($%
\textit{co}$N\smallsetminus N)\cap A=\emptyset ;$

\textit{iv) there exists }$M\in \langle A\rangle $\textit{\ such that}%
\newline
\textit{\ }$\bigcup\nolimits_{x\in M}[\bigcup\nolimits_{y\in S_{2}(x),\text{ 
}\exists t\in T(x,y)\text{ such that }F(t,y,x)\nsubseteq G(t,x,x)}($\textit{%
int}$_{X}\{u\in X:\exists t\in T(z_{x,y},u),$ \newline
$F(t,u,z_{x,y})\nsubseteq G(t,z_{x,y},z_{x,y})\}\cap S_{2}(x)]=X$\textit{;}

\textit{iv) for each }$u\in X,$\textit{\ the set }$M_{u}$\textit{\ is
convex, where}\newline
$M_{u}=\{x\in X:\exists t\in T(x,u),F(t,u,x)\nsubseteq G(t,x,x)\};$

\textit{v) for each }$x\in A,$\textit{\ }\newline
$\bigcup\nolimits_{\{y\in S_{2}(x),\text{ }t\in T(x,y)\text{: }%
F(t,y,x)\nsubseteq G(t,x,x)\}}($\textit{int}$_{X}\{u\in X:\exists t\in
T(z_{x,y},u),F(t,u,z_{x,y})\nsubseteq G(t,z_{x,y},z_{x,y})\}\cap
S_{2}(x)\subset \{u\in X:\exists t\in T(x,u),F(t,u,x)\nsubseteq G(t,x,x)\}$%
\textit{;}

\textit{Then, there exists }$x^{\ast }\in S_{1}(x^{\ast })$\ \textit{such
that }$\forall y\in S_{2}(x^{\ast }),$ $\forall t^{\ast }\in T(x^{\ast },y),$
$F(t^{\ast },y,x^{\ast })\subseteq G(t^{\ast },x^{\ast },x^{\ast })$.\textit{%
\medskip }

\textit{Proof.} Let $P,G:X\rightarrow 2^{X}$ be defined by

$P(x)=\{u\in X:\exists t\in T(x,u),F(t,u,x)\nsubseteq G(t,x,x)\},$ for each $%
x\in X$ and

$G(x)=S_{2}(x)\cap P(x),$ for each $x\in X.$

We are going to show that there exists $x^{\ast }\in X$ such that $x^{\ast
}\in S_{1}(x^{\ast })$ and $S_{2}(x^{\ast })\cap P(x^{\ast })=\emptyset .$

We consider two cases.

Case I.

$A=\{x\in X:$\textit{\ }$P(x)\cap S_{2}(x)\neq \emptyset \}=\{x\in X:$%
\textit{\ }$G(x)\neq \emptyset \}$\textit{\ }is nonempty.

The correspondence $G:A\rightarrow 2^{X},$ defined by $G(x)=S_{2}(x)\cap
P(x) $ for each $x\in A,$ is nonempty valued on $A.$

Further, let us define the correspondences $H,L,M:X\rightarrow 2^{X}$ by

$H(x)=\left\{ 
\begin{array}{c}
G(x),\text{ if }x\in A; \\ 
\emptyset ,\text{\ otherwise}%
\end{array}%
\right. ,$

$M(x)=\left\{ 
\begin{array}{c}
P(x),\text{ if }x\in A; \\ 
\emptyset ,\text{\ otherwise}%
\end{array}%
\right. $ and

$L(x)=\left\{ 
\begin{array}{c}
P(x),\text{ if }x\in A; \\ 
S_{1}(x),\text{\ otherwise.}%
\end{array}%
\right. $

For each $u\in X,$

$M^{-1}(u)=\{x\in X:$\textit{\ }$\mathit{u}\in M(x)\}=$

\ \ \ \ \ \ \ \ \ \ \ $=\{x\in A:$\textit{\ }$\mathit{u}\in P(x)\}=$

$\ \ \ \ \ \ \ \ =A\cap P^{-1}(u)=$

\ \ \ \ \ \ \ \ \ $=P^{-1}(u)=M_{u}.$

$L^{-1}(u)=\{x\in X:$\textit{\ }$\mathit{u}\in L(x)\}=$

\ \ \ \ \ \ \ \ \ \ \ $=\{x\in A:$\textit{\ }$\mathit{u}\in P(x)\}\cup
\{x\in X\setminus A:$\textit{\ }$\mathit{u}\in S_{1}(x)\}=$

$\ \ \ \ \ \ \ \ =(A\cap P^{-1}(u))\cup \lbrack (X\backslash A)\cap $\textit{%
\ }$(S_{1}^{-1}(u)]=$

\ \ \ \ \ \ \ \ \ $=P^{-1}(u)\cup \lbrack (X\backslash A)\cap
(S_{1}^{-1}(u)].$

According to assumption iv), $P^{-1}(u)$ is convex for each $u\in X$. Then,
co$M^{-1}(u)=$

=co$P^{-1}(u)=P^{-1}(u)\subset L^{-1}(u)$ for each $u\in X.$

According to i), $G$ is transfer open-valued. Assumptions i) and iii) imply
that there exists $M\in \langle A\rangle $ and for each $x\in M$ and $y\in
H(x),$ there exists $z_{x,y}\in X$ such that $y\in $int$_{X}H(z_{x,y})\cap
H(x)$ and $\tbigcup\nolimits_{x\in M}(\tbigcup\nolimits_{y\in H(x)}$int$%
_{X}H(z_{x,y}))=X$. In addition, assumption vi) implies $\tbigcup%
\nolimits_{y\in H(x)}$int$_{X}H(z_{x,y})\subseteq P(x)$ for each $x\in A.$

Let us define $Q:X\rightarrow 2^{X}$ by $Q(x):=X\backslash
\tbigcup\nolimits_{y\in H(x)}$int$_{X}H(z_{x,y})$ for each\textit{\ }$x\in
X. $

Then, $Q$ is closed-valued and $\tbigcap\nolimits_{x\in M}Q(x)=X\backslash
\tbigcup\nolimits_{x\in M}(\tbigcup\nolimits_{y\in H(x)}$int$%
_{X}H(z_{x,y}))=\emptyset ,$ where $M\subset A.$

According to Lemma 4, we conclude that $Q$ is not a KKM correspondence.
Thus, there exists $N\in \langle A\rangle $ such that\textit{\ }co$%
N\varsubsetneq Q(N)=\tbigcup\nolimits_{x\in N}(X\backslash
\tbigcup\nolimits_{y\in H(x)}$int$_{X}H(z_{x,y})).$

Hence, there exists $x^{\ast }\in $co$N$ with the property that $x^{\ast
}\in \tbigcup\nolimits_{y\in H(x)}$int$_{X}H(z_{x,y})\subset M(x)$ for each $%
x\in N.$ Therefore, there exists $x^{\ast }\in $co$N$ such that $x^{\ast
}\in M(x)$ for each $x\in N,$ which implies $N\subset M^{-1}(x^{\ast })$.
Further, it is true that co$N\subset $co$M^{-1}(x^{\ast })\subset
L^{-1}(x^{\ast })$. Consequently, $x^{\ast }\in $co$N\subset $co$%
M^{-1}(x^{\ast })\subset L^{-1}(x^{\ast }),$ which means that $x^{\ast }\in
L(x^{\ast }).$ We notice that, according to ii), $x\notin P(x)$ for each $%
x\in X,$ and then, $x^{\ast }\notin A.$ This fact is possible since $x^{\ast
}\in $co$N$ and assumption iii) asserts that $($co$N\smallsetminus N)\cap
A=\emptyset .$ Therefore, $x^{\ast }\in S_{1}(x^{\ast })$ and $G(x^{\ast
})=S_{2}(x^{\ast })\cap P(x^{\ast })=\emptyset .$

Consequently, there exists \textit{\ }$x^{\ast }\in S_{1}(x^{\ast })$\ such
that\textit{\ }$\forall y\in S_{2}(x^{\ast }),$ $\forall t^{\ast }\in
T(x^{\ast },y),$ $F(t^{\ast },y,x^{\ast })\subseteq G(t^{\ast },x^{\ast
},x^{\ast })$.

Case II.

$A=\{x\in X:$\textit{\ }$P(x)\cap S_{2}(x)\neq \emptyset \}=\{x\in X:$%
\textit{\ }$G(x)\neq \emptyset \}\mathit{=}\emptyset $.

In this case, $G(x)=\emptyset $ for each $x\in X.$ Let us define $%
Q:X\rightarrow 2^{X}$ by $Q(x):=X\backslash S_{1}(x)$ for each\textit{\ }$%
x\in X.$ The proof follows the same line as above and we obtain that there
exists $x^{\ast }\in X$ such that $x^{\ast }\in S_{1}(x^{\ast }).$

Obviously, $G(x^{\ast })=\emptyset .$ Consequently, the conclusion holds in
Case II.\smallskip

The next result can be obtained similarly as Theorem 23.

\begin{theorem}
\textit{Let }$Z$\textit{\ be a Hausdorff topological vector space and let }$%
X $\textit{\ be a nonempty compact convex subset of a\ topological vector
space }$E$\textit{. Let }$S_{1},S_{2}:X\rightarrow 2^{X},$\textit{\ }$%
T:X\times X\rightarrow 2^{X},$ $F:T(X\times X)\times X\times X\rightarrow
2^{Z}$\textit{\ and }$G:T(X\times X)\times X\times X\rightarrow 2^{Z}$ 
\textit{be correspondences with nonempty values. Assume that:}
\end{theorem}

\textit{i) }$S_{1},S_{2}$\textit{\ are open-valued and for each }$(x,y)\in
X\times X$\textit{\ with the property that }$\exists t\in T(x,y)$ such that $%
F(t,y,x)\cap G(t,x,x)\}=\emptyset ,$\textit{\ there exists }$z=z_{x,y}\in X$%
\textit{\ such that }$y\in $\textit{int}$_{X}\{u\in X:\exists t\in
T(z_{x,y},u),F(t,u,z_{x,y})\cap G(t,z_{x,y},z_{x,y})\}=\emptyset $\textit{; }

\textit{ii) for each }$x\in X$\textit{\ and }$t\in T(x,x),F(t,x,x)\subseteq
G(t,x,x);$\textit{\ }

iii) if $A=$ $\{x\in X:$\textit{\ there exist }$u\in S_{2}(x)$ \textit{and }$%
t\in T(x,u)$\textit{\ such that} $F(t,u,x)\cap G(t,x,x)=\emptyset \},$%
\textit{\ then for each }$N\in \langle A\rangle ,$\textit{\ }$($\textit{co}$%
N\smallsetminus N)\cap A=\emptyset ;$

\textit{iv) there exists }$M\in \langle A\rangle $\textit{\ such that}%
\newline
\textit{\ }$\bigcup\nolimits_{x\in M}[\bigcup\nolimits_{\{y\in S_{2}(x),%
\text{ }t\in T(x,y):\text{ }F(t,y,x)\nsubseteq G(t,x,x)\}}($\textit{int}$%
_{X}\{u\in X:\exists t\in T(z_{x,y},u),F(t,u,z_{x,y})\cap
G(t,z_{x,y},z_{x,y})=\emptyset \}\cap S_{2}(x)]=X$\textit{;}

\textit{iv) for each }$u\in X,$\textit{\ the set }$\{x\in X:\exists t\in
T(x,u),F(t,u,x)\cap G(t,x,x)=\emptyset \}$\textit{\ is convex}$;$

\textit{v) for each }$x\in A,$\newline
\textit{\ }$\bigcup\nolimits_{\{y\in S_{2}(x),\text{ }t\in T(x,y):\text{ }%
F(t,y,x)\nsubseteq G(t,x,x)\}}($\textit{int}$_{X}\{u\in X:\exists t\in
T(z_{x,y},u),F(t,u,z_{x,y})\cap G(t,z_{x,y},z_{x,y})=\emptyset \}\cap
S_{2}(x)\subset \{u\in X:\exists t\in T(x,u),F(t,u,x)\cap G(t,x,x)=\emptyset
\}$\textit{;}

\textit{Then, there exists }$x^{\ast }\in S_{1}(x^{\ast })$\ \textit{such
that }$\forall y\in S_{2}(x^{\ast }),$ $\forall t^{\ast }\in T(x^{\ast },y),$
$F(t^{\ast },y,x^{\ast })\cap G(t^{\ast },x^{\ast },x^{\ast })\neq \emptyset 
$.\textit{\medskip :}

\section{Concluding remarks}

In this paper, we have introduced $T$-properly quasi-convex correspondences
and $T$-properly quasi-convex sets. We have given some examples, as well. We
have used these notions to obtain coincidence-like theorems, to solve vector
equilibrium problems and to establish a generalized minimax inequality,
which is new in literature. We have also proved the existence of solutions
for vector variational inclusion problems, by applying the KKM principle.
Our research extends on some results which exist in literature. This study
can be continued by considering abstract convex spaces and generalized KKM
theorems.

\end{document}